\documentclass[%
a4paper,%
]%
{article}

\usepackage{amsmath}
\usepackage{amsthm}
\usepackage{amscd}

\usepackage{pxfonts}
\usepackage{xspace}

\newtheoremstyle{thm}{3pt}{3pt}{\itshape}{}{\bfseries}{}{.5em}{}
\newtheoremstyle{thmsub}{3pt}{3pt}{\upshape}{}{\bfseries}{}{.5em}{}
\theoremstyle{thm}
\newtheorem{theorem}{Theorem}[section]
\newtheorem{lemma}[theorem]{Lemma}
\newtheorem{proposition}[theorem]{Proposition}
\newtheorem{defn}[theorem]{Definition}
\newtheorem{corollary}[theorem]{Corollary}
\theoremstyle{thmsub}
\newtheorem{remark}[theorem]{Remark}

\newcommand{\m}[1]{\ensuremath {\mathcal{#1}}\xspace}
\newcommand{\mf}[1]{\ensuremath {\mathfrak{#1}}\xspace}
\newcommand{\ssetminus}{\ensuremath {\!\smallsetminus\!}}
\newcommand{\abs}[1]{\left\lvert #1 \right\rvert}
\newcommand{\norm}[1][\cdot]{\left\lVert #1 \right\rVert}
\newcommand{\ip}[2]{\langle #1, #2 \rangle}
\newcommand{\nmref}[1]{\ref{#1}}
\newcommand{\pdiff}[2]{\ensuremath\frac{\partial #1}{\partial #2}}
\newcommand{\cts}{{\text{cts}}}
\newcommand{\C}{\ensuremath {\mathbb{C}}\xspace}
\newcommand{\N}{\ensuremath {\mathbb{N}}\xspace}
\newcommand{\R}{\ensuremath {\mathbb{R}}\xspace}
\newcommand{\Ci}{\ensuremath {C^\infty}\xspace}
\newcommand{\exterior}{\Lambda}
\newcommand{\ad}{\text{ad}}
\newcommand{\ipv}[2]{\left(#1,#2\right)}

\DeclareMathOperator{\gl}{Gl}
\DeclareMathOperator{\cvx}{Cvx}
\DeclareMathOperator{\Diff}{Diff}
\DeclareMathOperator{\Ad}{Ad}
\DeclareMathOperator{\iso}{Iso}
\DeclareMathOperator{\map}{Map}

\numberwithin{equation}{section}

\begin{document}

\title{%
The Differential Topology of Loop Spaces
}  
\date{5th October 2005}
\author{%
Andrew Stacey%
}

\maketitle

\begin{abstract}
  This is an introduction to the subject of the differential topology
  of the space of smooth loops in a finite dimensional manifold.  It
  began as background notes to a series of seminars given at NTNU and
  subsequently at Sheffield.  The topics covered are: the smooth
  structure of the space of smooth loops; constructions involving
  vector bundles; submanifolds and tubular neighbourhoods; and a short
  introduction to the geometry and semi-infinite structure of loop
  spaces.  It is meant to be readable by anyone with a good grounding
  in finite dimensional differential topology.
\end{abstract}

\tableofcontents

\newpage

\section{Introduction}
\label{sec:intro}

This document started out as an accompaniment to a series of seminars
given at NTNU, and subsequently at Sheffield University, entitled:
``The Differential Topology of the Loop Space''.  The purpose of those
seminars was to present an introduction to this topic leading up to
the work contained in my preprint on the construction of a Dirac
operator on loop spaces, \cite{math.DG/0505077}.  This was intended to
be accessible to anyone with knowledge of basic finite dimensional
differential topology.  Thus there was a reasonable amount of
background material to be explained before the subject of Dirac
operators was broached.  The overall outline of this background
material was divided as follows:
\begin{enumerate}
\item The differential topology of loop spaces.

\item Spinors in arbitrary dimension.

\item The Dirac operator and the Atiyah-Singer index theorem in finite
  dimensions.
\end{enumerate}

The latter two topics are already superbly covered by books accessible
to any differential topologist, perhaps with a little functional
analysis.  The book~\cite{hlmm} is an excellent introduction to both
topics in finite dimensions whilst~\cite{rppr} deals with spinors in
arbitrary dimension.  A group-centric viewpoint is presented
in~\cite{apgs}, which is required reading for anyone seriously
thinking about loop spaces.

The main reference for the first topic is~\cite{akpm}.  Whilst several
of the standard texts on differential topology and geometry deal with
infinite dimensional manifolds, they only do so rigorously with
manifolds modelled on Banach spaces, which for loop spaces this
usually means either continuous or \(H^1\)-Sobolev loops.  Often
smooth and piece-wise smooth loops are used as a source of ``useful''
loops within these manifolds, but the space of smooth loops is not
given an actual manifold structure.  This is due to the technical
difficulties in extending calculus outside the realm of Banach spaces.

The purpose of~\cite{akpm} is to deal with precisely this issue and to
set up a theory of analysis in infinite dimensions in arbitrary
topological vector spaces.  It goes on to develop the theory of
infinite dimensional manifolds in its broadest setting which includes
smooth loop spaces.  However, this means that whilst being an
excellent book, its subject matter is perhaps too broad and too deep
for someone who just wants to know about the differential topology of
loop spaces.  The statement that the smooth loop space is a smooth
manifold appears in section~42, about two-thirds of the way through
the book.  Whilst the preceding four hundred pages are not all
required reading to get to this point, it is still a somewhat daunting
task to extract what one needs for loop spaces.

Therefore I started writing this document to fill in the details of
what I did not have time to say in my seminars.  The intention was to
provide a more gentle introduction than~\cite{akpm} but still include
the required detail to understand the differential topology of the
space of smooth loops.  It subsequently grew beyond that remit as I
thought of more topics that could be included without increasing the
technical difficulty.  A brief description of the topics covered is as
follows:

\begin{enumerate}
\item[\S \nmref{sec:smooth}] The smooth structure of the space of
  smooth loops in a finite dimensional manifold.

  We start by discussing what it means to be smooth outside the realm
  of Banach spaces, focusing particularly on the model spaces for
  loop spaces.  Using this we exhibit an atlas for the space of smooth
  loops and show that the transition functions are smooth.  Also in
  this section we show that various maps involving loop spaces are
  smooth.

\item[\S \nmref{sec:bundles}] The basic theory of vector bundles and
  associated principal and gauge bundles on the loop space defined by
  looping vector bundles on the original manifold.

  The main theme of this section is that looping something on the
  original manifold gives a corresponding object on the loop space.
  Thus the loop space of the tangent space is (diffeomorphic to) the
  tangent space of the loop space; vector bundles and their frame
  bundles loop to vector bundles and their frame bundles (in a sense)
  as do connections.  We conclude this section with an important
  example where this loopy behaviour fails: the loop space of the
  cotangent bundle is \emph{not} the cotangent bundle of the loop
  space.

\item[\S \nmref{sec:subtube}] Some important submanifolds with tubular
  neighbourhoods.

  We show that various submanifolds of a loop space have tubular
  neighbourhoods.  In particular, we consider submanifolds defined by
  imposing some condition on values that the loops can take at certain
  times -- for example, the based loop space -- and we consider fixed
  point submanifolds coming from the circle action.  One consequence
  of this is that the fundamental fibration \(\Omega M \to L M \to M\)
  is locally trivial.  We conclude with a submanifold that does not
  have a tubular neighbourhood.

\item[\S \nmref{sec:misc}] Basic introductions to two topics of further
  interest: geometry on loop spaces and the semi-infinite nature of
  loop spaces.

  We begin with a discussion of some of the issues in infinite
  dimensional geometry and prove some simple results concerning the
  Levi-Civita connection and geodesics.  We also give a very basic
  introduction to semi-infinite theory.
\end{enumerate}

This document began life as notes from talks given at NTNU and at
Sheffield so I would like to thank the topologists at those
institutions, and in particular Nils Baas, for letting me talk about
my favourite mathematical subject.  I would also like to thank Ralph
Cohen and the ``loop group'' at Stanford.

This is by no means a finished document, as an example it is somewhat
sparse on references.  Any comments, suggestions, and constructive
criticism will be welcomed.

\newpage

\section{Loopy Notation}

We first establish a little notation and some basic results for maps
between loop spaces that have nothing to do with calculus.  For a
(finite dimensional smooth) manifold \(M\), we denote \(\Ci(S^1,M)\)
by \(L M\).  If \(M\) is a space over \(S^1\), so that we have a
smooth map \(\pi : M \to S^1\), we denote the subspace of \(L M\)
consisting of sections of \(\pi\) by \(\Gamma_{S^1}(M)\).

\begin{defn}
  Let \(f : M \to N\) be a smooth map between manifolds.  Define \(f^L
  : L M \to L N\) by \(f^L(\alpha) = f \circ \alpha\).
\end{defn}

If \(M\) is a space over \(S^1\), \(f^L\) restricts to a map
\(\Gamma_{S^1}(M) \to L N\).  If, in addition, \(N\) is also a space
over \(S^1\) and \(f\) is a map of spaces over \(S^1\) then \(f^L\)
restricts to a map \(\Gamma_{S^1}(M) \to \Gamma_{S^1}(N)\).  We shall
use the notation \(f^L\) for all three of these maps.  Our convention
will be that, unless otherwise stated, \(f^L\) refers to the most
restrictive function that is applicable.

The most obvious type of space over \(S^1\) is the product: \(S^1
\times M\).  In this case, the projection \(p : S^1 \times M \to M\)
defines a bijection \(p^L : \Gamma_{S^1}(S^1 \times M) \to L M\).

\begin{defn}
  Let \(f : X \to L M\) be a map from a set \(X\) to the loop space of
  a manifold \(M\).  Let \(f^\lor : S^1 \times X \to M\) be the
  \emph{adjoint} of \(f\): \(f^\lor(t,x) = f(x)(t)\).
\end{defn}

\begin{defn}
  Let \(M\) be a smooth manifold.  Define the \emph{evaluation} map,
  \(e : S^1 \times L M \to M\) by \(e(t,\gamma) = \gamma(t)\).  We
  also write \(e_t : L M \to M\) for the \emph{evaluation at \(t\)}
  map, \(e_t(\gamma) = \gamma(t)\).
\end{defn}

\begin{lemma}
  \label{lem:loopadjprop}
  \begin{enumerate}
  \item The adjoint of \(f : X \to L M\) is the composition:
    \[
    f^\lor : S^1 \times X \xrightarrow{1 \times f} S^1 \times L M
    \xrightarrow{e} M.
    \]

  \item The adjoint of the identity on \(L M\) is the evaluation map.

  \item Let \(f : M \to N\) be a smooth map.  The adjoint of \(f^L\)
    is the map:
    \[
    S^1 \times L M \xrightarrow{e} M \xrightarrow{f} N
    \]

  \item Looping respects composition (i.e.~is a functor) in that \((f
    \circ g)^L = f^L \circ g^L\).  It also maps products to products
    in that \((f \times g)^L = f^L \times g^L\).

    \label{it:loopcomp}

  \item Let \(f : M \to N\), \(g : X \to L M\), and \(h : Y \to X\) be
    maps with \(f\) smooth.  The adjoint of \(f^L \circ g \circ h\) is
    \(f \circ g^\lor \circ (1 \times h)\).
  \end{enumerate}
\end{lemma}

\begin{proof}
  \begin{enumerate}
  \item We compute the composition as:
    \[
    (t,x) \xrightarrow{1 \times f} (t, f(x)) \xrightarrow{e} f(x)(t) =
    f^\lor(t,x).
    \]

  \item The follows from the above.

  \item We have:
    \[
    (f^L)^\lor(t, \alpha) = f^L(\alpha)(t) = (f \circ
    \alpha)(t) = f(\alpha(t)) = f \circ e(t, \alpha).
    \]
    Another way of writing this identity is: \(e \circ (1 \times f^L)
    = f \circ e\).

  \item This is obvious.

  \item Using the above, the adjoint of \(f^L \circ g \circ h\)
    is:
    \begin{align*}
      e \circ (1 \times (f^L \circ g \circ h)) &= e \circ (1 \times
      f^L)\circ (1 \times g) \circ (1 \times h) \\
      &= f \circ e\circ (1
      \times g) \circ (1 \times h) \\
      &= f \circ g^\lor \circ (1 \times
      h). \qedhere
    \end{align*}
  \end{enumerate}
\end{proof}

\newpage

\section{Smoothly Does It}
\label{sec:smooth}

The aim of this section is to prove that the space of smooth loops in
a finite dimensional manifold is again a smooth manifold.  After doing
this we go on to consider the topology of the space of smooth loops.
It may seem slightly topsy-turvey to consider the smooth structure
before the topology.  The reason is that in infinite dimensions the
link between continuity and smoothness is not as tight as it is in
finite dimensions (there are smooth maps which are not continuous).
After considering the topology we show that various obvious maps
involving loop spaces are smooth.  We conclude this section with a
short discussion of how to go about changing the type of loop under
consideration.

\subsection{Curvaceous Calculus}

The definitive reference for this section is the weighty
tome~\cite{akpm}.  Whilst a rewarding read for someone keen to learn
about infinite dimensional analysis, it is rather intricate for
someone just wanting to know the basics.  Therefore, we have tried in
this section to keep the material reasonably self-contained and
focused solely on the case of loop spaces.

The maxim in infinite dimensional analysis is:

\begin{quotation}
Smooth is as smooth does.
\end{quotation}

The problem is that once one strays outside the realm of Banach
spaces then there are many different ways of defining the derivative.
The \emph{Historical Remarks} at the end of \cite[ch I]{akpm} is a
good place to start investigating.  For example, one reference is
quoted as listing no less than \emph{twenty-five} inequivalent
definitions of the first derivative.  Even for Fr\'echet spaces -- the
next ``nicest'' spaces after Banach spaces, of which the space
\(\Ci(S^1, \R)\) is an example -- there are three inequivalent
definitions of infinite differentiability.  Therefore one needs to
consider what one wants from calculus before one actually defines it.

For example, smoothness and continuity are no longer as tightly linked
as they are in finite dimensions.  Let \(E\) be an infinite
dimensional locally convex topological vector space that is not
normable.  That is to say, the the topology on E cannot be given by a
single norm.  Consider the evaluation map \(E' \times E \to \R\),
\((f,e) \to f(e)\).  It is a straightforward result that with the
product topology on \(E' \times E\) then this map is not continuous.
However, one feels that it ought to be included in any ``reasonable''
calculus on \(E\).

Even in the theory of Hilbert spaces one finds that the topology on a
topological vector space is not fixed.  This becomes particularly
obvious when dealing with duality: on the dual space one can discuss
the strong topology, the weak topology, the simple topology, the
compact topology, and many, many more.  The guiding rule is to select
the topology most appropriate for the job in hand.

Thus with open season declared on the topology and a feeling of
freedom as to the definition of the derivative, one is left wondering
what exactly is well-defined.  The answer is: smooth curves.

It is simple to define smooth curves in any locally convex topological
vector space: a curve \(c : \R \to E\) is \emph{differentiable} if,
for all \(t\), the derivative \(c'(t)\) exists, where:
\[
c'(t) := \lim_{s \to 0} \frac1s \left( c(t + s) - c(t) \right).
\]
It is \emph{smooth} if all iterated derivatives exist.

Providing \(E\) satisfies a mild completeness condition (called
\emph{\(c^\infty\)-complete} or \emph{convenient} in \cite{akpm}),
then there is an incredibly simple recognition criterion for
smoothness of curves:
\begin{proposition}[{\cite[Theorem 2.14(4)]{akpm}}]
  Let \(E\) be a convenient locally convex space.  Then a curve \(c :
  \R \to E\) is smooth if and only if the curves \(l \circ c : \R \to
  \R\) are smooth for all \(l \in E^*\), the continuous dual of \(E\).
\end{proposition}

This appears to connect the concept of smoothness firmly to that of
continuity, but that is not so.  The smoothness, or otherwise, of a
curve depends only on what is called the \emph{bornology} of the space
\(E\).  That is, the family of \emph{bounded sets}.  The topology can
vary considerably without changing the bornology.

As this is our fixed point in analysis, it seems reasonable to define
everything else in terms of it.  This leads to:
\begin{defn}[{\cite[Definition 3.11]{akpm}}]
  A function \(f : E \supseteq U \to F\) defined on a
  \(c^\infty\)-open subset \(U\) of \(E\) is \emph{smooth} if it
  takes smooth curves in \(U\) to smooth curves in \(F\).
\end{defn}

I have sneaked in the term ``\(c^\infty\)-open'' about which I don't
intend to go into great detail.  Essentially, as we can vary the
topology somewhat without changing the notion of smoothness, it makes
sense to choose a particularly convenient topology.  This topology is
always at least as fine as the locally convex topology.  For Fr\'echet
spaces, it coincides with the standard topology.

To emphasise the point about smoothness being related to bornology
rather than topology, we record the following result about linear
maps:
\begin{lemma}[{\cite[Corollary 2.11]{akpm}}]
  A linear map \(l : E \to F\) between locally convex vector spaces is
  bounded, i.e.~bornological, if and only if it maps smooth curves in
  \(E\) to smooth curves in \(F\), i.e.~is a smooth map.
\end{lemma}
Thus as the evaluation map \(E' \times E \to \R\) is linear and
bounded it is smooth, but as mentioned above is not continuous for the
product topology.

The crucial tool in this calculus is the \emph{exponential law}.  One
important corollary of this is that in finite dimensions this notion
of smooth agrees with the standard one (originally proved in
\cite{jb3}).  To state the exponential law, we need first to note that
the space of smooth functions \(\Ci(U,F)\) has a natural topology
which is initial for the mappings \(c^* : \Ci(U,F) \to \Ci(\R,F)\) for
\(c \in \Ci(\R,U)\).

\begin{theorem}[{Exponential Law \cite[Theorem 3.12]{akpm}}]
  \label{th:exp}
  Let \(U_i \subseteq E_i\) be \(c^\infty\)-open subsets.  Then
  \(\Ci(U_1 \times U_2, F) \cong \Ci(U_1, \Ci(U_2, F))\).
\end{theorem}

This is called the exponential law because it can be written as
\(F^{U_1 \times U_2} \cong (F^{U_2})^{U_1}\).

We have now defined smoothness without reference to the derivative.
It is rather reassuring to discover that the derivative of a smooth
map is still definable:

\begin{theorem}[{\cite[Theorem 3.18]{akpm}}]
  Let \(E\) and \(F\) be locally convex spaces, and let \(U \subseteq
  E\) be \(c^\infty\)-open.  Then the differentiation operator
  \begin{gather*}
    d : \Ci(U,F) \to \Ci(U, \m{L}(E,F)), \\
    d f(x)(v) := \lim_{t \to 0} \frac{f(x + t v) - f(x)}{t},
  \end{gather*}
  exists and is bounded (smooth).  Also the chain rule
  holds:
  \[
  d( f \circ g)(x)v = d f(g(x)) d g(x) v.
  \]
\end{theorem}

Here, \(\m{L} (E, F)\) is the space of all \emph{bounded} linear maps
from \(E\) to \(F\).  We differ from \cite{akpm} in notation for this
space as we have reserved the letter ``\(L\)'' for loops and thus use
``\(\m{L}\)'' for linear maps.  It is a closed subspace of
\(\Ci(E,F)\) and is topologised as such.  Essentially, this result
says that the directional derivatives (which there is no problem
\emph{defining}, just problems with properties) fit together nicely to
define the global derivative.

\subsection{Euclidean Loop Spaces}

Before we can consider an arbitrary loop space we need to examine what
will become the model space, \(L \R^n\).  For those who like lists,
here is a list of the functional analysis properties that \(L \R^n\)
satisfies:
\begin{enumerate}
\item \(L \R^n\) is a Baire, barrelled, bornological, complete,
  convenient, Fr\'echet, infrabarrelled, Mackey, metrisable, Montel,
  nuclear, quasi-complete, reflexive, Schwartz, semi-reflexive, separable,
  space.  Also, the bidual map is a topological isomorphism, the
  \(c^\infty\) topology referred to above is the standard topology,
  and it has many nice properties involving tensor products.

\item the dual space, \((L \R^n)^*\) (with the strong topology), is
  all of the same except that it is not a Fr\'echet space -- and hence
  not metrisable -- nor is it a Baire space.  Also, the strong
  topology agrees with the Mackey topology and with the inductive
  topology.
\end{enumerate}

We are more interested in the calculus properties of this space,
although the ``nice properties involving tensor products'' are an
important factor in \cite{math.DG/0505077}.  We shall start with the
topology on \(L \R^n\).  As we noted in the previous section, this is
also the natural topology to consider when dealing with smooth maps on
\(L \R^n\).  To define it, we first define the topology on the space
\(C(S^1, \R^n)\) of continuous maps.

\begin{defn}
  A subbasis for the topology on \(C(S^1, \R^n)\) consists of the sets
  of the form:
  \[
  \{f \in C(S^1, \R^n) : f(K) \subseteq U\}
  \]
  where \(K \subseteq S^1\) is compact and \(U \subseteq \R^n\) is
  open.  This is called the \emph{compact-open} topology.
\end{defn}

In fact, \(C(S^1, \R^n)\) is a Banach space with supremum norm:
\[
\norm[\gamma]_\infty = \sup\{\norm[\gamma(t)] : t \in S^1\}.
\]

The topology on the smooth loop space is then defined as the topology
initial for the maps \(L \R^n \to C(S^1, \R^{k n})\) given by:
\[
\gamma \to \big(t \to (\gamma(t), \gamma'(t), \dotsc,
\gamma^{(k-1)}(t)) \big).
\]

If we replace continuous functions by some other type, for example:
square-integrable, we end up with the same topology on \(L \R^n\).
This is an example of how smoothness ``smooths out'' any initial
irregularities in definition.

We now turn to the smooth structure of \(L \R^n\).  Ultimately, we
need to know enough about the smooth maps on \(L \R^n\) to decide
whether or not the transition functions of an arbitrary loop space are
smooth.  These transition functions have a particularly simple
structure and so the task is easier than that of trying to
characterise \emph{all} smooth maps between open subsets of \(L
\R^n\).

We start by generalising the exponential law.  The problem with it as
stated in theorem~\ref{th:exp} is that the domains allowed are open
subsets of some linear space.  We wish to use \(S^1\) as one of the
domains.  That is, we wish to show:
\begin{proposition}
  \label{prop:exploop}
  A curve \(c : \R \to L \R^m\) is smooth if and only if its adjoint
  \(c^\lor : \R \times S^1 \to \R^m\) is smooth.
\end{proposition}

\begin{proof}
  Consider the quotient mapping \(\R \to S^1\).  This map completely
  determines the smooth structure of \(S^1\) in that a map \(S^1 \to
  \R^m\) is smooth if and only if the composite \(\R \to S^1 \to
  \R^m\) is smooth.  Similarly, \(\R \times S^1 \to \R^m\) is smooth
  if and only if \(\R \times \R \to \R^m\) is smooth.

  Thus given a curve \(c : \R \to L \R^m\) with adjoint \(c^\lor : \R
  \times S^1 \to \R^m\) we get a curve \(c_\R : \R \to \Ci(\R, \R^m)\)
  and a map \({c^\lor}_\R : \R \times \R \to \R^m\).  As both use the
  same projection \(\R \to S^1\), it is not hard to see that the
  adjoint of \(c_\R\) is \({c^\lor}_\R\); that is, \({c_\R}^\lor =
  {c^\lor}_\R\).  Thus by the standard exponential theorem,
  theorem~\ref{th:exp}, \(c_\R\) is smooth if and only if
  \({c^\lor}_\R\) is smooth.  By the above, \({c^\lor}_\R\) is smooth
  if and only if \(c^\lor\) is smooth.  Thus it remains to show that
  \(c\) is smooth if and only if \(c_\R\) is smooth.

  Now the map \(\R \to S^1\) exhibits \(L \R^m\) as a linear subspace
  of \(\Ci(\R, \R^m)\) in that the induced map \(L \R^m \to \Ci(\R,
  \R^m)\) is a homeomorphism onto its image.  The image is obviously a
  closed subspace and hence \(c^\infty\)--closed.  Therefore, by
  \cite[Lemma 3.8]{akpm}, a curve in \(L \R^m\) is smooth if and only
  if it is smooth as a curve in \(\Ci(\R, \R^m)\).
\end{proof}

The type of function that we shall examine is the following: let \(V
\subseteq S^1 \times \R^n\) and \(W \subseteq S^1 \times \R^m\) be
open subsets, and let \(\psi : V \to W\) be a map of spaces over
\(S^1\); here, \(S^1 \times \R^k\) has the obvious structure of a
space over \(S^1\) and a subspace of a space over \(S^1\) has the
obvious inherited structure.

\begin{lemma}
  \label{lem:smloop}
  Under the natural identification \(\Gamma_{S^1}(S^1 \times \R^k)
  \cong L \R^k\), \(\Gamma_{S^1}(V)\) and \(\Gamma_{S^1}(W)\) are open
  subsets of the respective loop spaces and \(\psi^L\) is a smooth map.
\end{lemma}

Note that \(\Gamma_{S^1}(V)\) is empty unless the map \(V \to S^1\) is
surjective.

\begin{proof}
  That \(\Gamma_{S^1}(V)\) identifies with an open subset of \(L
  \R^n\) is obvious as it is true for the topology induced by the
  inclusion into the space of continuous loops.

  To show that \(\psi^L\) is smooth, it is sufficient to take the case
  where \(W = S^1 \times \R^m\).  This means that we are considering
  the target to be \(L \R^n\) via the bijection \(p^L :
  \Gamma_{S^1}(S^1 \times \R^m) \to L \R^m\).  Thus we need to show
  that if \(c : \R \to \Gamma_{S^1}(V)\) is smooth then \(p^L \circ
  \psi^L \circ c : \R \to L \R^m\) is smooth.

  By the exponential law, it is sufficient to consider the adjoint of
  this.  Using lemma~\ref{lem:loopadjprop} we see that this is \(p
  \circ \psi \circ c^\lor : \R \times S^1 \to \R^m\).  This is
  smooth and hence \(p^L \circ \psi^L\) takes smooth maps to smooth
  maps.  Thus \(\psi^L\) is smooth.
\end{proof}

Let \(d_v \psi : V \times \R^n \to W \times \R^m\) be the vertical
derivative of \(\psi\).

\begin{lemma}
  \label{lem:deriv}
  \(d(\psi^L) = (d_v \psi)^L\).
\end{lemma}

\begin{proof}
  Recall that the derivative of a smooth function is determined by the
  directional derivatives.  Using the identification of
  \(\Gamma_{S^1}(S^1 \times \R^k)\) with \(L \R^k\), the derivative of
  \(\psi^L\) is a smooth map \(\Gamma_{S^1}(V) \to \m{L}(L \R^n, L
  \R^m)\).  So the derivative at a loop \(\alpha \in \Gamma_{S^1}(V)\)
  in the direction of \(\beta \in L \R^n\) is an element of \(L \R^m\)
  which is the limit of the following net indexed by \(s \in \R^+\):
  \[
  \frac{\psi^L(\alpha + s \beta) - \psi^L(\alpha)}{s}.
  \]
  Now evaluation at time \(t\) is a continuous linear map \(L \R^n \to
  \R^n\) which takes this net to:
  \[
  \frac{\psi^L(\alpha + s \beta)(t) - \psi^L(\alpha)(t)}{s} =
    \frac{\psi(t, \alpha(t) + s \beta(t)) - \psi(t, \alpha(t))}{s}.
  \]
  This is precisely the difference quotient which tends to the
  vertical derivative of \(\psi\) at \(\alpha(t)\) in the direction
  \(\beta(t)\).  Since a loop is completely determined by its values
  at each time,
  \(d(\psi^L)(\alpha)\beta\) is the same loop as \(t \to d_v \psi(t,
  \alpha(t))\beta(t)\).  Hence:
  \[
  d(\psi^L)(\alpha)\beta = (d_v \psi)^L(\alpha) \beta. \qedhere
  \]
\end{proof}

\begin{corollary}
  The map \(d(\psi^L)(\alpha) : L \R^n \to L \R^m\) is \(L
  \R\)--linear.
\end{corollary}

\begin{proof}
  Since \(d \psi^L(\alpha)\) is \R-linear, the part to prove in the
  statement about linearity is that \(d \psi^L(\alpha)\) commutes with
  the action of \(L \R\).  Let \(\nu \in L \R\), then:
  \begin{align*}
  d \psi^L(\alpha)(\nu \beta)(t) &= d_v \psi(t, \alpha(t))(\nu(t)
  \beta(t)) \\
  &= \nu(t) d_v \psi(t, \alpha(t))\beta(t) \\
  &= \big( \nu d \psi^L(\alpha) \beta \big)(t).
  \end{align*}
  Hence \(d \psi^L(\alpha)(\nu \beta) = \nu d \psi^L(\alpha) \beta\)
  as required.
\end{proof}

\begin{lemma}
  \label{lem:diffeom}
  Now suppose that \(n = m\) and that \(\psi : V \to W\) is a
  diffeomorphism with inverse \(\phi : W \to V\).  Then \(\psi^L :
  \Gamma_{S^1}(V) \to \Gamma_{S^1}(W)\) is a diffeomorphism with
  inverse \(\phi^L\).
\end{lemma}

\begin{proof}
  Both \(\psi^L\) and \(\phi^L\) are smooth.  By
  property~\ref{it:loopcomp} of lemma~\ref{lem:loopadjprop}, the
  compositions \(\psi^L \circ \phi^L\) and \(\phi^L \circ \psi^L\) are
  the identities on their domains.  Hence \(\psi^L\) is a
  diffeomorphism with inver \(\phi^L\).
\end{proof}

This has a few nice consequences that we shall exploit.  Most are
things that jolly well ought to be true, but given that the definition
of smooth is probably unfamiliar it's as well to spell out the
details.

\begin{corollary}
  \label{cor:vbdiff}
  \begin{enumerate}
  \item Let \(\alpha : S^1 \to \Diff(\R^n)\) be a smooth map (the
    group \(\Diff(\R^n)\) inherits a smooth structure from \(\Ci(\R^n,
    \R^n)\)).  Then the induced map \(L \R^n \to L \R^n\) given by
    \(\gamma \to (t \to \alpha(t)\gamma(t))\) is a diffeomorphism.

  \item Let \(E \to S^1\) be a smooth, orientable vector bundle.  Then
    the space of sections, \(\Gamma_{S^1}(E)\), has a natural smooth
    structure.

  \item The results above generalise to orientable smooth vector
    bundles as follows: let \(E, F \to S^1\) be orientable smooth
    vector bundles.  Let \(V \subseteq E\) and \(W \subseteq F\) be
    open subsets.  Let \(\psi : V \to W\) be a smooth map covering the
    identity on \(S^1\).  Then \(\psi^L : \Gamma_{S^1}(V) \to
    \Gamma_{S^1}(W)\) is smooth with derivative the loop of the
    vertical derivative of \(\psi\).  If \(\psi\) is a diffeomorphism
    with inverse \(\phi\) then \(\psi^L\) is a diffeomorphism with
    inverse \(\phi^L\).
  \end{enumerate}
\end{corollary}

The key in the second is to choose a diffeomorphism of (the total
space of) \(E\) with \(S^1 \times \R^n\) and so identify the space of
sections of \(E\) with \(L \R^n\).  This defines a smooth structure on
this space of sections which is independent of the diffeomorphism
chosen by the first result.

The same results hold for non-orientable bundles but are a little 
fiddly to prove.  One has to consider what are called \emph{twisted}
loop spaces.  These add no analytic difficulties, but perhaps add one
or two conceptual ones so we shall avoid using them.

\subsection{The Smooth Structure of an Arbitrary Loop Space}
\label{sec:loopstr}

We are now able to consider the smooth loop space of a finite
dimensional smooth manifold, \(M\).  We shall make two assumptions on
the type of manifold that we consider: one necessary and one for
convenience.
\begin{enumerate}
\item \(\partial M = \emptyset\).

  The loop space of a manifold with boundary is a complicated object.
  To put it simply, it is not true that \(L M = L (M \ssetminus
  \partial M) \cup L (\partial M)\).  The best description of \(L M\)
  in this case is as a stratified space, with strata indexed by closed
  subspaces of the circle.  The layer corresponding to \(F \subseteq
  S^1\) consists of those loops \(\alpha : S^1 \to M\) such that
  \(\alpha^{-1}(\partial M) = F\).  The top level, corresponding to
  the empty set, is \(L(M \ssetminus \partial M)\).  The next level
  consists of all loops which intersect \(\partial M\) at one point,
  and so on.

\item \(M\) is orientable.

  This allows us to stay out of the twisted realm.  We shall not
  actually use this in the analysis, it is merely to make this account
  reasonably self-contained.
\end{enumerate}

The key tool for defining the charts for the loop space is the notion
of a \emph{local addition} on \(M\), cf~\cite[\S 42.4]{akpm}:
\begin{defn}
\label{def:locadd}
  A \emph{local addition} on \(M\) consists of a smooth map \(\eta : T
  M \to M\) such that
  \begin{enumerate}
  \item the composition of \(\eta\) with the zero section is the
    identity on \(M\), and

  \item there exists an open neighbourhood \(V\) of the diagonal in
    \(M\) such that the map \(\pi \times \eta : T M \to M \times M\) is
    a diffeomorphism onto \(V\).
  \end{enumerate}
\end{defn}

In \cite[\S 42.4]{akpm}, the above is called a \emph{globally defined
local addition} but the difference is not important for us.  We shall
later want to relax this further by replacing the tangent bundle by an
arbitrary vector bundle (which must, \emph{a fortiori}, be isomorphic
to the tangent bundle) but for now we stick with the tangent bundle to
keep things conceptually simple.  The following result is contained in
the discussion following the definition of a local addition in
\cite[\S 42.4]{akpm}:

\begin{proposition}
  Any finite dimensional manifold without boundary admits a local
  addition.
\end{proposition}

\begin{proof}
  Without going into great detail, the essentials of the proof are
  that the exponential map coming from a Riemannian structure almost
  defines a local addition except that the domain of the
  diffeomorphism is not the whole tangent space (except in a few
  simple cases) but a neighbourhood of the zero section.  The proof is
  completed by exhibiting a smooth fibre-preserving embedding of the
  total space of the tangent bundle into the domain of the
  diffeomorphism.  The composition of this with the exponential map is
  the required local addition.
\end{proof}

Let \(\eta : T M \to M\) be a local addition on \(M\).  Let \(V
\subseteq M \times M\) be the image of the map \(\pi \times \eta : T M
\to M \times M\).  Although as yet we know nothing about the
topologies of \(L T M\) or of \(L V\), we can at least say that the
looped map, \((\pi \times \eta)^L\), is a bijection.

\begin{lemma}
  Let \(\alpha \in L M\).  Define the set \(U_\alpha \subseteq L M\)
  by:
  \[
  U_\alpha := \{\beta \in L M : (\alpha, \beta) \in L V\}.
  \]
  Then the preimage of \(\{\alpha\} \times U_\alpha\) under \((\pi
  \times \eta)^L\) is naturally identified with
  \(\Gamma_{S^1}(\alpha^* T M)\).  In particular, the zero section of
  \(\alpha^* T M\) maps to \((\alpha,\alpha) \in \{\alpha\} \times
  U_\alpha\).
\end{lemma}

\begin{proof}
We claim that there is a diagram:
\[
\begin{CD}
  L T M @> (\pi \times \eta)^L >> L V \\
  @AAA @AA \beta \to (\alpha,\beta) A \\
  \Gamma_{S^1}(\alpha^* T M) @. U_\alpha,
\end{CD}
\]
such that the bijection at the top takes the image of the left-hand
vertical map to the image of the right-hand one.  Both of the vertical
maps are injective -- the right-hand one obviously so, we shall
investigate the left-hand one in a moment -- and thus the bijection
\((\pi \times \eta)^L\) induces a bijection from the lower left to the
lower right.

The left-hand vertical map, \(\Gamma_{S^1}(\alpha^* T M) \to L T M\),
is defined as follows: the total space \(\alpha^* T M\) is:
\[
\{(t,v) \in S^1 \times T M : \alpha(t) = \pi(v)\}.
\]
It is an embedded submanifold of \(S^1 \times T M\).  Therefore a map
into \(\alpha^* T M\) is smooth if and only if the compositions with
the projections to \(S^1\) and to \(T M\) are smooth.  Now a map \(S^1
\to \alpha^* T M\) is a section if and only if it projects to the
identity on \(S^1\).  Therefore there is a bijection (of sets):
\begin{align*}
\Gamma_{S^1}(\alpha^* T M) &\cong \{\beta \in L T M : (t,
\beta(t)) \in \alpha^* T M \text{ for all } t \in S^1\} \\
&= \{\beta \in L T M : \alpha(t) = \pi \beta(t) \text{ for all
} t \in S^1\} \\
&= \{\beta \in L T M : \pi^L \beta = \alpha\} =: L_\alpha T M.
\end{align*}
In particular, the map \(\Gamma_{S^1}(\alpha^* T M) \to L T M\) is
injective.

We apply \((\pi \times \eta)^L\) to the defining condition for
\(L_\alpha T M\) and see that \(L_\alpha T M\) is the preimage under
this map of everything of the form \((\alpha, \gamma)\) in \(L V\).
By construction, \(\gamma \in L M\) is such that \((\alpha, \gamma)
\in L V\) if and only if \(\gamma \in U_\alpha\).  Hence \((\pi \times
\eta)^L\) identifies \(L_\alpha T M\) with \(\{\alpha\} \times
U_\alpha\).

Finally, note that the zero section of \(\alpha^* T M\) maps to the
image of \(\alpha\) under the zero section of \(T M\).  Since \(\eta\)
composed with the zero section of \(T M\) is the identity on \(M\),
the image of the zero section of \(\alpha^* T M\) in \(V\) is
\((\alpha,\alpha)\) as required.
\end{proof}

\begin{defn}
  Let \(\Psi_\alpha : \Gamma_{S^1}(\alpha^* T M) \to U_\alpha\) be the
  resulting bijection.
\end{defn}

In detail, this map is as follows: let \(\beta \in
\Gamma_{S^1}(\alpha^* T M)\) and let \(\tilde{\beta}\) be the
corresponding loop in \(T M\), so \(\beta(t) = (t,
\tilde{\beta}(t))\).  Then \((\pi \times \eta)^L(\tilde{\beta}) =
(\alpha, \eta^L(\tilde{\beta}))\) so \(\Psi_\alpha(\beta) =
\eta^L(\tilde{\beta})\).

The domains of these functions are sections over \(S^1\) of smooth
orientable finite dimensional vector bundles.  We equip these with the
smooth structure as in corollary~\ref{cor:vbdiff}.  On the other side,
since \(\alpha \in U_\alpha\), the family \(\{U_\alpha\}\) covers the
loop space of \(M\).  Thus we have an open covering of the loop space
by sets each of which we can identify with a linear space with a
chosen smooth structure.

We now turn to the investigation of the transition functions.  We
shall prove a slightly stronger result than is needed at this moment
in that we shall allow the local addition to vary.  As well as showing
that the transition functions are smooth this will show that the
maximal atlas so defined contains all such charts with any choice of
local addition.

Let \(\eta_1, \eta_2 : T M \to M\) be local additions with
corresponding neighbourhoods \(V_1, V_2\) of the diagonal in \(M
\times M\).  Let \(\alpha_1, \alpha_2\) be smooth loops in \(M\).  Let
\(\Psi_1 : \Gamma_{S^1}(\alpha_1^* T M) \to U_1\), \(\Psi_2 :
\Gamma_{S^1}(\alpha_2^* T M) \to U_2\) be the corresponding charts.

\begin{lemma}
  Let \(W_{1 2} \subseteq \alpha_1^* T M\) be the set:
  \[
  \{(t,v) \in \alpha_1^* T M : (\alpha_2(t), \eta_1(v)) \in V_2\}.
  \]
  Then \(W_{1 2}\) is open and \({\Psi_1}^{-1}(U_1 \cap U_2) =
  \Gamma_{S^1}(W_{1 2})\).
\end{lemma}

\begin{proof}
  The set \(W_{1 2}\) is open as it is the preimage of an open set via
  a continuous map.  To show the second statement we need to prove
  that \(\gamma \in \Gamma_{S^1}(\alpha_1^* T M)\) takes values in
  \(W_{1 2}\) if and only if \(\Psi_1(\gamma) \in U_2\) (by
  construction we already have \(\Psi_1(\gamma) \in U_1\)).

  So let \(\gamma \in \Gamma_{S^1}(\alpha^* T M)\) and  let
  \(\tilde{\gamma}\) be the image of \(\gamma\) in \(L T M\).
  Thus \(\gamma(t) = (t, \tilde{\gamma}(t))\).
  Now \(\gamma\) takes values in \(W_{1 2}\) if and only if:
  \[
  \big(\alpha_2(t), \eta_1(\tilde{\gamma}(t))\big) \in V_2
  \]
  for all \(t \in S^1\).  That is to say, if and only if \((\alpha_2,
  {\eta_1}^L(\tilde{\gamma})) \in L V_2\).  By definition, this is
  equivalent to the statement that \({\eta_1}^L(\tilde{\gamma}) \in
  U_2\).  Now \({\eta_1}^L(\tilde{\gamma}) = \Psi_1(\gamma)\) so
  \(\gamma\) takes values in \(W_{1 2}\) if and only if
  \(\Psi_1(\gamma) \in U_1 \cap U_2\).
\end{proof}

\begin{proposition}
  \label{prop:trans}
  The transition function:
  \[
  \Phi_{1 2} := {\Psi_1}^{-1} \Psi_2 : {\Psi_1}^{-1}(U_{1 2}) \to
  {\Psi_2}^{-1}(U_{1 2})
  \]
  is a diffeomorphism.
\end{proposition}

\begin{proof}
  We define \(W_{2 1} \subseteq \alpha_2^* T M\) as the set of \((t,v)
  \in \alpha_2^* T M\) with \((\alpha_1(t), \eta_2(v)) \in V_1\).  As
  for \(W_{1 2}\), \(\Psi_2^{-1}(U_1 \cap U_2) = \Gamma_{S^1}(W_{2
  1})\).

  The idea of the proof is to set up a diffeomorphism between \(W_{1
  2}\) and \(W_{2 1}\).  Corollary~\ref{cor:vbdiff} says that the
  resulting map on sections is a diffeomorphism.  Finally, we show
  that this diffeomorphism is the transition function defined in the
  statement of this proposition.

  Let \(\theta_1 : W_{1 2} \to T M\) be the map:
  \[
  \theta_1(t,v) = (\pi \times \eta_2)^{-1}(\alpha_2(t), \eta_1(v)).
  \]
  The definition of \(W_{1 2}\) ensures that \((\alpha_2(t),
  \eta_1(v)) \in V_2\) for \((t,v) \in W_{1 2}\) and this is the image
  of \(\pi \times \eta_2\).  Hence \(\theta_1\) is well-defined.  Define
  \(\theta_2 : W_{2 1} \to T M\) similarly.  These are both smooth maps.

  Notice that \(\pi(\pi \times \eta_i)^{-1} : V_i \subseteq M \times
  M \to M\) is the projection onto the first factor and \(\eta_i(\pi
  \times \eta_i)^{-1} : V_i \to M\) is the projection onto the second.
  Thus \(\pi \theta_1(t,v) = \alpha_2(t)\).  Hence \(\theta_1 : W_{1
  2} \to T M\) is such that \((t, \theta_1(t,v)) \in \alpha_2^* T M\)
  for all \((t,v) \in W_{1 2}\).  Then:
  \[
  \big(\alpha_1(t), \eta_2(\theta_1(t,v))\big) = (\alpha_1(t), \eta_1(v))
  \in V_1
  \]
  so \((t, \theta_1(t,v)) \in W_{2 1}\).  Hence we have a map
  \(\phi_{1 2} : W_{1 2} \to W_{2 1}\) given by:
  \[
  \phi_{1 2}(t,v) = (t, \theta_1(t,v)).
  \]
  Similarly we have a map \(\phi_{2 1} : W_{2 1} \to W_{1 2}\).  These
  are both smooth since the composition with the inclusion into \(S^1
  \times T M\) is smooth.

  Consider the composition \(\phi_{2 1}\phi_{1 2}(t,v)\).  Expanding
  this out yields:
  \begin{align*}
    \phi_{2 1}\phi_{1 2}(t,v) &= \phi_{2 1}(t, \theta_1(t,v)) \\
    &= (t, \theta_2(t,\theta_1(t,v))) \\
    &= \big(t, (\pi \times \eta_1)^{-1}(\alpha_1(t),
    \eta_2(\theta_1(t,v))) \big) \\
    &= \big(t, (\pi \times \eta_1)^{-1}(\alpha_1(t), \eta_1(v)) \big)
    \\
    &= \big(t, (\pi \times \eta_1)^{-1}(\pi(v), \eta_1(v)) \big) \\
    &= (t, v).
  \end{align*}
  The penultimate line is because \((t,v) \in \alpha_1^* T M\) so
  \(\pi(v) = \alpha_1(t)\).

  Hence \(\phi_{2 1}\) is the inverse of \(\phi_{1 2}\) and so
  \(\phi_{1 2} : W_{1 2} \to W_{2 1}\) is a diffeomorphism.  Thus the
  map \({\phi_{1 2}}^L\) is a diffeomorphism from \({\Psi_1}^{-1}(U_1
  \cap U_2)\) to \({\Psi_2}^{-1}(U_1 \cap U_2)\).  We just need to
  show that this is the transition function.  It is sufficient to show
  that \(\Psi_2 {\phi_{1 2}}^L = \Psi_2 \Phi_{1 2}\).  The right-hand
  side is, by definition, \(\Psi_1\), which satisfies:
  \[
  \Psi_1(\gamma)(t) = \eta_1(\tilde{\gamma}(t))
  \]
  where \(\tilde{\gamma} : S^1 \to T M\) is such that \(\gamma(t)
  = (t, \tilde{\gamma}(t))\).  On the other side:
  \begin{align*}
    {\phi_{1 2}}^L(\gamma)(t) &= \phi_{1 2}(\gamma(t)) \\
  &= \big(t, \theta_1(t, \tilde{\gamma}(t))\big) \\
  &= \big(t, (\pi \times \eta_2)^{-1}(\alpha_2(t),
  \eta_1(\tilde{\gamma}(t))) \big),
  \intertext{hence:}
  \Psi_2 {\phi_{1 2}}^L(\gamma)(t) &= \eta_2(\pi \times
  \eta_2)^{-1}( \alpha_2(t), \eta_1(\tilde{\gamma}(t))) \\
  &= \eta_1(\tilde{\gamma})(t),
  \end{align*}
  as required.  Thus \({\phi_{1 2}}^L = \Phi_{1 2}\) and so
  the transition functions are diffeomorphisms.
\end{proof}

\begin{remark}
  There is no \emph{a priori} reason why we needed to assume that
  \(M\) was a finite dimensional manifold.  All we needed was to know
  what the smooth structure of \(M\) was and to know that it had a
  local addition.  Since the loop of a local addition is again a local
  addition, we could iterate this construction and show that any
  iterated loop space is a smooth manifold.
\end{remark}

\subsection{The Loop Space of the Tangent Space}

The next topic that we wish to consider is the topology of the space
\(L M\).  It is usually defined, by analogy with the topology of \(L
\R^n\), as a projective limit of the maps \(L M \to C(S^1, T^{(k)}
M)\) where \(T^{(k)} M\) is recursively defined as the tangent space
of \(T^{(k-1)} M\), and \(T^{(1)} M = T M\).  We shall not define the
topology this way but shall show that our definition is equivalent.
To do this, we need to know how the structure of \(L M\) relates to
that of \(L T M\).  This is a useful piece of knowledge that we shall
use again when we look at the tangent space of \(L M\) so we record it
in its own section here.

We mentioned earlier that we wish to replace the source of the local
addition by an arbitrary vector bundle, albeit one isomorphic to the
tangent bundle.  We can either repeat all of the above with arbitrary
sources for the two local additions to see that the charts so defined
lie in the same atlas, or we can simply observe that as the source
must be isomorphic to the tangent bundle, we can use an isomorphism to
transfer the local addition to the tangent bundle and also to define a
linear diffeomorphism between the model spaces.  Thus any chart
defined using an arbitrary vector bundle factors through one that uses
the tangent bundle and hence is in the same atlas.

The first application of this is to lift the charts defined for \(L
M\) to charts for \(L T M\).  Let \(T^{(2)} M\) be the tangent bundle
of \(T M\).  Let \(\eta : T M \to M\) be a local addition with
corresponding open neighbourhood \(V \subseteq M \times M\).  Then \(d
\eta : T^{(2)} M \to T M\) is a local addition in this new sense but
\emph{not} in the strict sense.

This is because \(d \pi \times d \eta : T^{(2)} M \to T V\) is a
diffeomorphism and \(T V\) is an open neighbourhood in \(T M \times T
M\) of the diagonal.  However while \(d \pi : T^{(2)} M \to T M\) is a
vector bundle projection it is not the projection of a tangent bundle
onto its base.  Rather we have a commutative diagram of vector bundle
projections:
\[
\begin{CD}
  T^{(2)} M @>d\pi>> T M \\
  @V\pi VV @VV\pi V \\
  T M @>\pi>> M
\end{CD}
\]
where every map \emph{except} \(d\pi\) is the natural projection of a
tangent bundle onto its base.

In all of the following, we consider \(T^{(2)} M\) as a vector bundle
over \(T M\) via the map \(d \pi\).

Using this local addition we get a chart \(\Phi_\alpha :
\Gamma_{S^1}(\alpha^* T^{(2)} M) \to V_\alpha\) for any \(\alpha \in L
T M\).  The first thing to observe about this chart is that \(V_\alpha
= (\pi^L)^{-1}(U_{\pi^L \alpha})\), where \(\pi^L : L T M \to L M\) is
the loop of the projection and \(U_{\pi^L \alpha}\) is the codomain of
the chart at \(\pi^L \alpha\) defined using the local addition
\(\eta\).  This is because \(T V = (\pi \times \pi)^{-1} V\) and so:
\begin{align*}
  V_\alpha &= \{\beta \in L T M : (\alpha(t), \beta(t)) \in T V \text{
  all } t \in S^1\} \\
  &= \{\beta \in L T M : (\pi \alpha(t), \pi \beta(t)) \in V \text{
  all } t \in S^1\} \\
  &= \{\beta \in L T M : \pi^L \beta \in U_{\pi^L\alpha}\} \\
  &= (\pi^L)^{-1} U_{\pi^L \alpha}.
\end{align*}
Thus to get a cover of \(L T M\), we start with a family of loops in
\(M\) such that the corresponding charts cover \(L M\) and then we
pick any-old lifts of these loops to \(T M\).  One of the most
straightforward ways to choose such a lift is as follows: let \(\tau :
S^1 \to T S^1\) be a section.  Define \(\hat{\tau} : L M \to L T M\)
by \(\hat{\tau}\alpha = d\alpha \circ \tau : S^1 \to T S^1 \to T M\).

\begin{proposition}
  \label{prop:charttm}
  The map \(\hat{\tau}\) induces a natural map on sections:
  \(\hat{\tau} : \Gamma_{S^1}(\alpha^* T M) \to
  \Gamma_{S^1}((\hat{\tau} \alpha)^* T^{(2)} M)\).  A trivialisation
  of \(\alpha^* T M\) induces a trivialisation of \((\hat{\tau}
  \alpha)^* T^{(2)} M\) such that the following diagram commutes:
  \[
  \begin{CD}
    L \R^n @>>> \Gamma_{S^1}(\alpha^* T M) @>\Psi_\alpha>> U_\alpha
    @>\subseteq>> L M \\
    @V\hat{\tau}VV @V\hat{\tau}VV @. @V\hat{\tau}VV \\
    L T\R^n @>>> \Gamma_{S^1}((\hat{\tau}\alpha)^* T^{(2)} M)
    @>\Phi_{\hat{\tau}\alpha}>> V_{\hat{\tau}\alpha} @>\subseteq>> L T
    M.
  \end{CD}
  \]
\end{proposition}

\begin{proof}
  We shall start by explaining the two maps induced by \(\tau\).  The
  space \(\alpha^* T M\) is an embedded submanifold of \(S^1 \times T
  M\).  It can be defined as the preimage of the diagonal in \(M
  \times M\) under the map \((t,v) \to (\alpha(t), \pi(v))\).  Its
  tangent space is therefore an embedded submanifold of \(T S^1 \times
  T^{(2)} M\) and is the preimage of the diagonal in \(T M \times T
  M\) under the map \(d \alpha \times d \pi\).  Hence \(T(\alpha^* T
  M)\) is the pull-back of \(d \pi : T^{(2)} M \to T M\) under \(d
  \alpha\).

  Now a section \(\beta\) of \(\alpha^* T M\) is a special kind of map
  \(\beta : S^1 \to \alpha^* T M\).  It therefore differentiates to a
  map \(d \beta : T S^1 \to T(\alpha^* T M) = (d\alpha)^* T^{(2)} M\).
  As \(\beta\) is a section, \(\pi \beta\) is the identity on \(S^1\),
  hence \(d \pi d \beta\) is the identity on \(T S^1\) and so \(d
  \beta\) is a section of \((d\alpha)^* T^{(2)} M\).  We compose this
  with \(\tau\) to get \(d \beta \circ \tau : S^1 \to (d\alpha)^*
  T^{(2)} M\).  As \(d \beta\) is a section, \(d \pi d \beta \tau =
  \tau\) and so \(d \beta \tau\) is a section of \(\tau^* (d\alpha)^*
  T^{(2)} M\).  As \(\tau^* (d \alpha)^* = (d \alpha \tau)^*\), \(d
  \beta \tau\) is a section of of \((\hat{\tau}\alpha)^* T^{(2)} M\).
  Thus we have the map \(\hat{\tau} : \beta \to d\beta \tau\) from
  sections of \(\alpha^* T M\) to sections of \((\hat{\tau}\alpha)^*
  T^{(2)} M\).

  The other map works similarly: let \(\phi : S^1 \times \R^n \to
  \alpha^* T M\) be a trivialisation.  This differentiates to \(d \phi
  : T S^1 \times T \R^n \to T(\alpha^* T M) = (d \alpha)^* T^{(2)} M\)
  which is a trivialisation over \(T S^1\).  This pulls-back via
  \(\tau\) to a trivialisation \(S^1 \times T \R^n \to \tau^* (d
  \alpha)^* T^{(2)} M\).  The left-hand side of this is \(S^1 \times
  \R^{2 n}\) and the right is \((\hat{\tau}\alpha)^* T^{(2)} M\).

  To see how the commutative diagram works, we start with a loop in
  \(\R^n\) and trace how it becomes a loop in \(M\).  Let
  \(\check{\beta} : S^1 \to S^1 \times \R^n\) be the map \(t \to (t,
  \beta(t))\), then the loop in \(M\) is the composition:
  \[
  S^1 \xrightarrow{\check{\beta}} S^1 \times \R^n \xrightarrow{\phi}
  \alpha^* T M \subseteq S^1 \times T M \to T M \xrightarrow{\eta} M.
  \]

  We differentiate this, and use the identification of \(T(\alpha^* T
  M)\) with \((d \alpha)^* T^{(2)} M \subseteq T S^1 \times T^{(2)}
  M\) to get the map:
  \[
  T S^1 \xrightarrow{d\check{\beta}} T S^1 \times T \R^n
  \xrightarrow{d \phi} (d\alpha)^* T M \subseteq T S^1 \times T^{(2)}
  M \to T^{(2)} M \xrightarrow{d \eta} T M.
  \]

  The first four terms in this sequence are fibre bundles over \(T
  S^1\) (viewing \(T S^1\) as a bundle over itself with one-point
  fibres) and all the maps are of fibre bundles.  Therefore we can
  pull-back these terms over \(S^1\) using the section \(\tau : S^1
  \to T S^1\).  The map from the fourth to fifth terms is the
  projection of the product away from \(T S^1\) so is unchanged under
  the pull-back.  The third term is \((\hat{\tau}\alpha)^* T^{(2)} M\)
  and the second map is the induced trivialisation of this bundle.
  The first map, \(\tau^* d \check{\beta}\), fits into the diagram:
  \[
  \begin{CD}
    T S^1 @> d\check{\beta} >> T S^1 \times T \R^n \\
    @A\tau AA @A \tau \times 1 AA \\
    S^1 @> \tau^* d \check{\beta} >> S^1 \times T \R^n.
  \end{CD}
  \]
  Differentiating \(\check{\beta}(t) = (t,\beta(t))\) yields
  \(d\check{\beta} = (1, d\beta)\) and so \(d \check{\beta} \circ \tau
  = (\tau, d \beta \circ \tau) = (\tau, \hat{\tau} \beta)\).  Hence
  \((\tau^* d \beta) (t) = (t,\hat{\tau} \beta(t))\).  Thus \(\tau^* d
  \beta\) is the map \(S^1 \to T \R^n\) corresponding to the map
  \(\hat{\tau} \beta \in L T \R^n\).  Write this as \(\hat{\tau}
  \check{\beta}\).

  We therefore have the commutative diagram:
  \[
  \begin{CD}
    T S^1 @>d \check{\beta}>> T S^1 \times T \R^n @>d\phi>> (d
    \alpha)^* T^{(2)} M @>>> T^{(2)} M @> d\eta >> T M \\
    @A\hat{\tau}AA @A\hat{\tau}AA @A\hat{\tau}AA @A1AA
    @A1AA \\
    S^1 @>\hat{\tau}\check{\beta}>> S^1 \times T \R^n @>\tau^* d\phi>>
    (\hat{\tau} \alpha)^* T^{(2)} M @>>> T^{(2)} M @> d \eta >> T M.
  \end{CD}
  \]
  Going first up and then all the way along results in \((d
  \Psi_\alpha(\beta)) \tau = \hat{\tau} \Psi_\alpha(\beta)\).  Going
  up at the third stage, \((\hat{\tau} \alpha)^* T^{(2)} M\), shows
  that this factors through the map \(\hat{\tau}\) on sections, whilst
  going all the way along and then up shows that it factors through
  the map \(\hat{\tau} : L \R^n \to L T \R^n\).  Hence the commutative
  diagram of the statement of this proposition is true.
\end{proof}

\subsection{The Topology of the Loop Space}  

We now turn to the topology on \(L M\).

\begin{defn}
  The topology on \(L M\) is the topology such that \(W \subseteq L
  M\) is open if and only if \(\Psi_\alpha^{-1}(W)\) is open for all
  \(\Psi_\alpha\).
\end{defn}

Since the transition functions are diffeomorphisms and hence
homeomorphisms, it is sufficient to check this condition using a
family of charts which cover \(L M\).

To investigate this topology we start with the topology of the
continuous loop space, \(C(S^1, M)\).  The topology on this is similar
to the ``compact-open'' topology on \(C(S^1, \R^n)\).  A basis
consists of the sets:
\[
N(K,U) := \{\alpha \in C(S^1, M) : \alpha(K) \subseteq U\}
\]
where \(K\) runs through the family of compact subsets of \(S^1\) and
\(U\) through the open subsets of \(M\).

With this topology, the space \(C(S^1, M)\) has the following
properties:
\begin{enumerate}
\item It is a separable metrisable space.  An explicit metric is
  \[
  d_\cts(\alpha, \beta) := \sup\{d_M(\alpha(t), \beta(t)) : t \in
  S^1\}
  \]
  where \(d_M\) is a metric on \(M\) compatible with its topology; for
  example, that which comes from a Riemannian structure.

\item The space of smooth loops is dense.  The standard proof of this
  is to show that any continuous loop can be approximated to an
  arbitrary degree by a piecewise geodesic.  Any piecewise smooth loop
  can be approximated (in the compact-open topology) by a smooth one
  by ``rounding off the corners''.
\end{enumerate}

\begin{lemma}
  The natural inclusion \(L M \to C(S^1, M)\) is continuous and the
  codomains of the charts are open for the induced topology.
\end{lemma}

\begin{proof}
  We start with the statement on the codomains.  Recall that the
  codomain of \(\Psi_\alpha\) is the set \(U_\alpha\) such that
  \(\{\alpha\} \times U_\alpha = (\{\alpha\} \times L M) \cap L V\).
  Now a loop in \(V\) is smooth if and only if it is smooth as a loop
  in \(M \times M\), so as \(\alpha\) is smooth the right-hand side of
  this is:
  \[
  \Big(\big(\{\alpha\} \times C(S^1, M) \big) \cap C(S^1, V) \Big)
  \cap (L M \times L M).
  \]
  Hence if \(U^0_\alpha \subseteq C(S^1, M)\) is the continuous
  version of \(U_\alpha\), then \(U_\alpha = L M \cap U^0_\alpha\).
  As \(V\) is open in \(M \times M\), \(C(S^1, V)\) is open in
  \(C(S^1, M)\) and so \(U^0_\alpha\) is open in \(C(S^1, M)\).  Hence
  \(U_\alpha\) is open for the induced topology on \(L M\).

  The topology on \(L M\) is defined such that \(W \subseteq L M\) is
  open if and only if \(\Psi_\alpha^{-1}(W)\) is open for all
  \(\alpha\).  Thus a map \(f\) from \(L M\) is continuous if and only
  if all the compositions \(f \circ \Psi_\alpha\) are continuous.

  A smooth trivialisation of \(\alpha^* T M\) identifies
  \(\Gamma_{S^1}(\alpha^* T M)\) with \(L \R^n\).  It also identifies
  the space of continuous sections, \(\Gamma_{S^1}^0(\alpha^* T M)\),
  with the continuous loop space, \(C(S^1, \R^n)\).  The obvious
  square commutes.  As with smooth sections, a section of \(\alpha^* T
  M\) is continuous if and only if the induced map \(S^1 \to T M\) is
  continuous.  Therefore there is an analogous map:
  \[
  \Psi_\alpha^0 : \Gamma_{S^1}^0(\alpha^* T M) \to U^0_\alpha, \qquad
  \Psi_\alpha^0(\beta)(t) = \eta(\widetilde{\beta}(t)).
  \]
  In this case, we have a topology on the target already.  That the
  map \(\Psi_\alpha^0\) is continuous for this topology is a direct
  consequence of the fact that if \(f : X \to Y\) is continuous then
  the induced map \(f^L : C(S^1, X) \to C(S^1, Y)\) is continuous.

  Consider the following diagram:
  \[
  \begin{CD}
    L \R^n @>\cong >> \Gamma_{S^1}(\gamma^* T M) @>\Psi_\gamma >> L M
    \\
    @Vi_{R^n}VV @Vi_\gamma VV @Vi_MVV \\
    C(S^1, \R^n) @> \cong >> \Gamma_{S^1}^0 (\gamma^* T M) @>
    \Psi_\gamma^0 >> C(S^1, M)
  \end{CD}
  \]

  The map \(i_M\) is continuous if and only if \(i_M \circ
  \Psi_\gamma\) is continuous, for each \(\gamma\).  This is the same
  as \(\Psi^0_\gamma i_\gamma\).  Then \(i_\gamma\) is continuous
  because it is equivalent to \(i_{\R^n}\) which is the inclusion of
  smooth loops in continuous loops in Euclidean space.
\end{proof}

This result has some useful corollaries:

\begin{corollary}
  With this topology, \(L M\) is:
  \begin{enumerate}
  \item Hausdorff,
  \item regular,
  \item second countable, and
  \item paracompact.
  \end{enumerate}
\end{corollary}

\begin{proof}
  \begin{enumerate}
  \item Any topology finer than a Hausdorff topology is Hausdorff.
    Viewing \(L M\) as a (topological) subspace of \(C(S^1,M)\), the
    induced topology is Hausdorff as \(C(S^1, M)\) is metrisable.
    Hence the given topology on \(L M\) is Hausdorff.

  \item With the compact-open topology, \(L M\) is a regular space.
    Therefore, as the codomains of charts are open for the
    compact-open topology, every point has a closed neighbourhood
    contained in the codomain of a chart.

    Let \(C \subseteq L M\) be a closed set and \(\gamma \notin C\).
    Let \(D\) be a closed neighbourhood of \(\gamma\) that is
    contained in the codomain of a chart.  The chart map induces a
    homeomorphism between \(D\) and a closed subset of a locally
    convex topological vector space, which is regular as a topological
    space.

    Therefore, there is a set \(\widetilde{B}\) which is open in \(D\)
    such that \(\gamma \in \widetilde{B}\) and the closure of
    \(\widetilde{B}\) does not meet \(D \cap C\).  Let \(B\) be the
    intersection of \(\widetilde{B}\) with the interior of \(D\).  As
    \(D\) is a neighbourhood of \(\gamma\), this is still a non-empty
    set containing \(\gamma\).  Moreover, as an open subset of the
    interior of \(D\), it is open in \(L M\).  Its closure in \(L M\)
    is the same as its closure in \(D\).  This does not meet \(D \cap
    C\), whence -- as it is contained within \(D\) -- it does not meet
    \(C\).

    Hence \(L M\) is regular.

  \item We start by observing that a countable number of codomains of
    charts will cover \(L M\).  This is because with the compact-open
    topology, i.e.~with \(L M\) viewed as a subspace of \(C(S^1, M)\),
    \(L M\) is second countable (as it is a subspace of a second
    countable space).  It is therefore Lindel\"of.  Now the codomains
    of the charts form an open covering \emph{with the compact-open
    topology}, hence a countable number of them, say the family
    \(\{U_n\}\), will cover \(L M\).

    Each chart codomain is homeomorphic to an open subset of the
    separable, metrisable space \(L \R^n\).  It therefore is itself
    second countable.  Thus for each codomain in our countable
    covering we choose a countable basis for the topology.  This forms
    a countable family of open sets which we claim is a basis for the
    topology.

    This is straightforward: for \(W \subseteq L M\) open, as the
    countable family of codomains covers \(L M\) then \(W\) is the
    union of \(W \cap U_n\).  Each of these is then the union of sets
    in our basis.  Hence \(W\) is the union of sets in our basis.
    Thus \(L M\) has a countable basis and so is second countable.

  \item We use the following two results:
    \begin{enumerate}
    \item A second countable space is Lindel\"of, and
    \item A Lindel\"of regular space is paracompact. \qedhere
    \end{enumerate}
  \end{enumerate}
\end{proof}

From this, various topological results follow including the fact that
\(L M\) is metrisable.  We shall quote from \cite{akpm} the result
that \(L M\) is \emph{smoothly paracompact} which means that it admits
smooth partitions of unity subordinate to any open cover.  This
implies that it is \emph{smoothly Hausdorff} in that any two points
can be separated by a smooth function.

We have now proved that:

\begin{theorem}
  The space \(L M\) is a smooth manifold.
\end{theorem}

The topology on \(L M\) is usually given in a different way to that
which we have used.  The section \(\tau : S^1 \to T S^1\), \(\tau(t) =
(t, \pdiff{}{t})\) defines a section \(L M \to L T M\) and thus, by
iteration, \(L M \to L T^{(k)} M\).  Composing this with the inclusion
of smooth loops in continuous loops yields a family of maps \(L M \to
C(S^1, T^{(k)} M)\).

\begin{proposition}
  \label{prop:projtop}
  The topology on \(L M\) is the projective topology for this family
  of maps \(L M \to C(S^1, T^{(k)} M)\).
\end{proposition}

\begin{proof}
  We apply proposition~\ref{prop:charttm} with the section \(\tau(t) =
  (t, \pdiff{}{t})\) to see that the following diagram commutes:
  \[
  \begin{CD}
    L \R^n @>\Psi_\alpha>> L M \\
    @V\beta \to \beta'VV @VV \gamma \to \gamma' V \\
    L T \R^n @>\Phi_\alpha>> L T M \\
    @VVV @VVV \\
    C(S^1, T \R^n) @>>> C(S^1, T M).
  \end{CD}
  \]
  The topology on \(L \R^n\) is the projective limit of the obvious
  iteration of this diagram.  Therefore we can construct a basis for
  the topology on \(L M\) out of sets which are the preimage of open
  sets in each \(C(S^1, T^{(k)} M)\), which is precisely what is meant
  by the definition of the projective topology.
\end{proof}

\subsection{Loopy Maps}

In this section we look at some maps involving loop spaces that ought
to be smooth and show that indeed they are so.

One corollary of the proof of proposition~\ref{prop:projtop} is that
the canonical map \(L M \to L T M\) given by \(\alpha \to \alpha'\) is
continuous.  Further, we can show that it is smooth:

\begin{lemma}
  Let \(\tau : S^1 \to T S^1\) be a section.  The map \(L M \to L T
  M\), \(\alpha \to \hat{\tau} \alpha = (d \alpha) \tau \), is smooth.
\end{lemma}

\begin{proof}
  Let \(\alpha \in L M\).  From the proof of
  proposition~\ref{prop:projtop} we have charts for \(L M\) and \(L T
  M\) with chart maps \(\Psi_\alpha\) and \(\Phi_{\hat{\tau}\alpha}\)
  such that the following diagram commutes:
  \[
  \begin{CD}
    L \R^n @>\Psi_\alpha>> L M \\
    @V{\beta \to \hat{\tau}\beta}VV @VV{\gamma \to \hat{\tau}\gamma}V
    \\
    L T \R^n @>\Phi_{\hat{\tau}\alpha} >> L T M.
  \end{CD}
  \]
  Hence if we can show that the statement of the lemma is true for \(M
  = \R^n\) we can deduce that it holds for all \(M\).  However, this
  is obvious as it is a linear bounded map.
\end{proof}

There is another class of maps that it will be useful to know are
smooth:

\begin{proposition}
  \label{prop:loopmap}
  Let \(f : M \to N\) be a smooth map between smooth finite
  dimensional manifolds then \(f^L : L M \to L N\) is smooth.
\end{proposition}

\begin{proof}
  This is essentially the same argument as when we showed that the
  transition functions were smooth.  We need to show that the loop of
  \(f\) is smooth in charts.  Let \(\eta_M : T M \to M\) and \(\eta_N
  : T N \to N\) be local additions with neighbourhoods \(V_M\) and
  \(V_N\) of the respective diagonals.  Let \(\alpha \in L M\).  We
  therefore have charts \(\Psi_M : \alpha^* T M \to U_M\) and \(\Psi_N
  : (f^L\alpha)^* T N \to U_N\) at \(\alpha\) and \((f^L\alpha)\)
  respectively.

  Let \(W \subseteq \alpha^* T M\) be the set:
  \[
  \{(t,v) \in \alpha^* T M : (f(\alpha(t)), f(\eta_M(v))) \in V_N\}.
  \]
  We claim that \(\beta \in \Gamma_{S^1}( \alpha^* T M)\) takes values
  in \(W\) if and only if \(f^L( \Psi_M(\beta)) \in U_N\).  Let
  \(\tilde{\beta} \in L T M\) be such that the section \(t \to (t,
  \tilde{\beta}(t))\) is \(\beta\).  Then \(\Psi_M(\beta) =
  {\eta_M}^L(\tilde{\beta})\).  By the definition of \(U_N\),
  \(f^L(\Psi_M(\beta)) \in U_N\) if and only if \((f(\alpha)(t),
  f(\Psi_M(\beta))(t)) \in V_N\) for all \(t\).  As \(f(\alpha)(t) =
  f(\alpha(t))\) and \(f(\Psi_M(\beta))(t) =
  f(\eta_M(\tilde{\beta}(t)))\), this is equivalent to \((t,
  \tilde{\beta}(t)) \in W\) for all \(t\), i.e.~that \(\beta\) takes
  values in \(W\).

  As for the transition functions, let \(\theta : W \to T N\) be the
  smooth function:
  \[
  \theta(t,v) = (\pi_N \times \eta_N)^{-1} (f(\alpha(t)),
  f(\eta_M(v))).
  \]
  Exactly as for the transition functions, \(\pi_N \theta(t,v) =
  f(\alpha(t))\) and so the map \((t,v) \to (t,\theta(t,v))\) is a
  smooth map \(W \to (f^L\alpha)^* T N\) (smooth as it is smooth into
  \(S^1 \times T N\)).

  Hence we have a smooth map \(\Theta\) from sections that take values
  in \(W\) into \(\Gamma_{S^1}(f(\alpha)^* T N)\).  It remains to show
  that this is the map induced by \(f\); that is, \(f^L \Psi_M = \Psi_N
  \Theta\).  This is straightforward.  Let \(\beta \in
  \Gamma_{S^1}(\alpha^* T M)\) take values in \(W\) and let
  \(\tilde{\beta}\) be the related loop in \(T M\).  Then:
  \[
  (f^L \Psi_M)(\beta)(t) = f (\eta_M( \tilde{\beta}(t)))
  \]
  whilst:
  \[
  (\Psi_N \Theta)(\beta) = \eta_N(\pi_N \times \eta_N)^{-1}
  (f(\alpha(t)), f(\eta_M(\tilde{\beta}(t)))) =
  f(\eta_M(\tilde{\beta}(t))).
  \]
  Hence \(\Theta\) is the map on charts corresponding to \(f^L\).  Thus
  \(f^L : L M \to L N\) is smooth.
\end{proof}

The other result about maps that we wish to prove is a variant of the
exponential law:

\begin{theorem}
  \label{th:expmfd}
  Let \(M, N\) be smooth manifolds with \(M\) finite dimensional.  Let
  \(f : N \to L M\) be a map and let \(f^\lor : S^1 \times N \to M\)
  be the adjoint: \(f^\lor(t,x) = f(x)(t)\).  Then \(f\) is smooth
  if and only if \(f^\lor\) is smooth.  That is, the assignment \(f
  \to f^\lor\) is an identification:
  \[
  \Ci(N, L M) \cong \Ci(S^1 \times N, M).
  \]
\end{theorem}

As stated, this is a corollary of \cite[Theorem 42.14]{akpm} with
their \((\m{M}, M, \m{N})\) corresponding to our \((N, S^1, M)\).  We
give a self-contained proof here.

\begin{proof}
  We start with the linear case: \(M = \R^n\) and \(N\) is an open
  subset of a convenient vector space.  Thus we have \(f : N \to L
  \R^n\) and \(f^\lor : S^1 \times N \to \R^n\).

  Let \(\pi : \R \to S^1\) be the standard covering map.  Now \(S^1
  \times N\) is a smooth manifold modelled on \(\R \times N\) using
  the usual charts for \(S^1\) so \(\pi \times 1 : \R \times N \to S^1
  \times N\) is a local diffeomorphism.  Hence \(f^\lor\) is smooth if
  and only if \(f^\lor(\pi \times 1)\) is smooth.

  On the other side, pre-composition with \(\pi\) defines a linear
  embedding of \(L \R^n\) onto a closed subspace of \(\Ci(\R, \R^n)\)
  as in the proof of proposition~\ref{prop:exploop}.  Therefore \(f\)
  is smooth if and only if the map \(x \to f(x) \circ \pi\) is smooth.

  The adjoint of \(x \to f(x) \circ \pi\) is the map:
  \[
  (t,x) \to (f(x) \circ \pi)(t) = f(x)(\pi(t)) = f^\lor(\pi(t),x) =
  f^\lor(\pi \times 1)(t,x).
  \]
  Thus by the original exponential law, theorem~\ref{th:exp}, the map
  \(x \to f(x) \circ \pi\) is smooth if and only if \(f^\lor(\pi
  \times 1)\) is smooth.  Hence \(f\) is smooth if and only if
  \(f^\lor\) is smooth.

  In the general case, we use charts to move into the linear one.
  However, with things as they stand we find that we are only
  partially successful: we may assume that \(N\) is an open subset of
  a convenient vector space and that the codomain of \(f\) is
  contained in the codomain of a chart but it may not be possible to
  assume that the image of \(f^\lor : S^1 \times N \to M\) is
  contained in the codomain of a chart.  The easiest way around this
  is to enhance \(f^\lor\) to the map \(f^\land : S^1 \times N \to
  S^1 \times M\), \((t,x) \to (t, f^\lor(t,x))\).  This is smooth if
  and only if \(f^\lor\) is smooth.

  Let \(\Psi_\alpha : \Gamma_{S^1}(\alpha^* T M) \to U_\alpha\) be a
  chart for \(L M\) defined using a local addition, \(\eta\), on \(M\)
  with neighbourhood \(V\) of the diagonal in \(M \times M\).  Let
  \(V_\alpha \subseteq S^1 \times M\) be the open set \(\{(t,x) :
  (\alpha(t), x) \in V\}\).  We claim that the evaluation \(S^1 \times
  U_\alpha \to S^1 \times M\) maps into \(V_\alpha\) and moreover that
  there is a commutative diagram:
  \begin{equation}
    \label{diag:eval}
  \begin{CD}
    S^1 \times \Gamma_{S^1}(\alpha^* T M) @>1 \times \Psi_\alpha>> S^1
    \times U_\alpha \\
    @VVV @VVV \\
    \alpha^* T M @>\psi_\alpha>> V_\alpha
  \end{CD}
  \end{equation}
  where both vertical maps are evaluations and both horizontal maps
  are diffeomorphisms.  The map \(\psi_\alpha\) is defined as follows:
  there is a map
  \[
  \alpha^* T M \subseteq S^1 \times T M \xrightarrow{1 \times \eta} S^1
  \times M
  \]
  which, we claim, has image \(V_\alpha\) and is a diffeomorphism onto
  that image.  To see this, consider the diagram:
  \[
  \begin{CD}
    S^1 \times T M @>1 \times \pi \times \eta>> S^1 \times V
    @>\subseteq >> S^1 \times M \times M \\
    @AAA @AAA @AAA \\
    \alpha^* T M @. V_\alpha @>\subseteq>> S^1 \times M,
  \end{CD}
  \]
  where the right-hand vertical map is \((t,x) \to (t,\alpha(t),x)\).
  Each vertical map is an embedding and the first upper horizontal map
  is a diffeomorphism.  It is clear that \(\alpha^* T M\) corresponds
  to \(V_\alpha\) on the upper level and hence the map \(\alpha^* T M
  \to V_\alpha\) is a diffeomorphism.

  To see that diagram~\eqref{diag:eval} commutes, let \(\beta \in
  \Gamma_{S^1}(\alpha^* T M)\) and \(\tilde{\beta} \in L T M\) be
  related as usual, so \(\Psi_\alpha(\beta) = \eta(\tilde{\beta})\).
  Thus \(\Psi_\alpha(\beta)(t) = \eta(\tilde{\beta})(t) =
  \eta(\tilde{\beta}(t))\) and so going along and then down is the
  map:
  \[
  (t, \beta) \to (t, \eta(\tilde{\beta}(t))).
  \]
  Going down and then along, we get:
  \[
  (1 \times \eta)(\beta(t)) = (1 \times \eta)(t, \tilde{\beta}(t)) =
  (t, \eta(\tilde{\beta}(t))).
  \]

  Now \(f^\land\) is the composition of:
  \[
  S^1 \times N \xrightarrow{1 \times f} S^1 \times U_\alpha \to
  V_\alpha \subseteq S^1 \times M
  \]
  where the map \(S^1 \times U_\alpha \to V_\alpha\) is the map
  \((t,\gamma) \to (t, \gamma(t))\).  Thus the pair \((f, f^\land)\)
  defines a pair \((g, g^\land)\) with \(g : N \to
  \Gamma_{S^1}(\alpha^* T M)\) and \(g^\land : S^1 \times N \to
  \alpha^* T M\) such that \(g\) is smooth if and only if \(f\) is
  smooth and \(g^\land\) smooth if and only if \(f^\land\) is
  smooth.

  A smooth trivialisation of \(\alpha^* T M\) identifies that space
  with \(S^1 \times \R^n\) and sections with \(L \R^n\).  Thus we can
  regard \(g\) and \(g^\land\) as maps into \(L \R^n\) and \(S^1
  \times \R^n\) respectively.  Moreover, \(g^\land(t,x) = (t,
  g^\lor(t,x))\) so \(g^\land\) is smooth if and only if \(g^\lor :
  S^1 \times N \to \R^n\) is smooth.

  We are thus in the linear case and so \(g\) is smooth if and only if
  \(g^\lor\) is smooth.  Hence \(f\) is smooth if and only if
  \(f^\lor\) is smooth.
\end{proof}

\begin{corollary}
  Let \(e : S^1 \times L M \to M\) be the evaluation map:
  \(e(t,\gamma) = \gamma(t)\).  Let \(\iota : M \to L M\) be the
  inclusion of constant loops: \(\iota(x) = (t \to x)\).  Both \(e\)
  and \(\iota\) are smooth.
\end{corollary}

\begin{proof}
  The map \(e\) is adjoint to the identity \(L M \to L M\), hence is
  smooth, whereas \(\iota\) adjoints to the projection map \(S^1
  \times M \to M\), hence is smooth.
\end{proof}

For the second of those maps we can actually prove a stronger result:

\begin{proposition}
  \label{prop:constembed}
  The map \(\iota : M \to L M\) is an embedding.
\end{proposition}

\begin{proof}
  Let \(\eta : T M \to M\) be a local addition on \(M\) with \(V
  \subseteq M \times M\) the corresponding neighbourhood of the
  diagonal.  For \(x \in M\) let \(\eta_x : T_x M \to V_x \subseteq
  M\) be the restriction of \(\eta\) to \(T_x M\), where \(V_x =
  \eta(T_x M)\) is such that \(\{x\} \times V_x = (\{x\} \times M)
  \cap V\).  It is a simple result from finite dimensional
  differential topology that this family defines an atlas for \(M\).

  Let \(\Psi_x : \Gamma_{S^1}({\gamma_x}^* T M) \to U_x\) be the chart
  for \(L M\) at the constant loop \(\gamma_x\).  As \(\gamma_x\) is a
  constant loop, \({\gamma_x}^* T M = S^1 \times T_x M\).  Define
  \(T_x M \to \Gamma_{S^1}({\gamma_x}^* T M)\) by \(v \to \beta_v\)
  where \(\beta_v(t) = (t,v)\).  This is smooth since it is linear and
  bounded.

  Now \(\{x\} \times V_x = (\{x\} \times M) \cap V\) and
  \(\{\gamma_x\} \times U_x = (\{\gamma_x\} \times L M) \cap L V\).
  Under the inclusion \(M \to L M\) and \(V \to L V\) we thus have
  \(V_x = M \cap U_x\).  We claim that the following diagram commutes:
  \[
  \begin{CD}
    T_x M @>\eta_x>> V_x \\
    @VVV @VVV \\
    \Gamma_{S^1}({\gamma_x}^* T M) @>\Psi_x >> U_x.
  \end{CD}
  \]
  To see this, let \(\beta \in \Gamma_{S^1}({\gamma_x}^* T M)\) and
  \(\tilde{\beta} \in L T M\) be related as usual.  Thus
  \(\tilde{\beta}\) satisfies \(\pi \tilde{\beta} = \gamma_x\) and so
  \(\tilde{\beta} \in L T_x M\).  By definition, \(\Psi_x(\beta) =
  \eta(\tilde{\beta}) = \eta_x(\tilde{\beta})\).  Now
  \(\tilde{\beta}_v\) is the constant map at \(v \in T_x M\) and so
  \(\Psi_x(\beta_v)\) is the constant map at \(\eta_x(v)\).  Hence the
  above diagram commutes.
  
  Thus these charts satisfy the requirements for exhibiting \(M\) as a
  submanifold of \(L M\).
\end{proof}

\subsection{Vector Space Idol}

In proving the topological results about the space of smooth loops we
used the fact that the chart maps extend naturally to continuous
loops.  In fact, the chart maps can be used to construct a manifold
structure on the spaces of many different types of loop.  The key
ingredients are:
\begin{enumerate}
\item The type of loop is at least continuous (with topology at least
  as fine as the compact-open topology) and at most smooth, and

\item The smooth structure is diffeomorphism-invariant.  That is to
  say, for \(V, W \subseteq S^1 \times \R^n\) open subsets and
  \(\psi : V \to W\) a diffeomorphism, the induced map \(\psi^L\) on
  the appropriate type of sections is a diffeomorphism.
\end{enumerate}

With these two properties, we can adapt the work of
section~\ref{sec:loopstr} to this new type of loop.  To keep the
alterations as straightforward as possible, we still anchor our charts
at genuinely smooth loops: the first condition above ensures that
these charts still cover the loop space.  If we allowed other anchors,
we would have to expand the second condition above to homeomorphisms
that are fibrewise diffeomorphisms and satisfy some sort of condition
amongst the fibres.

Where we may run into trouble is with the topological properties of
the new loop space.  In particular, we need the model space to be
second countable to prove that the loop space is second countable.  If
this fails, paracompactness is also thrown into doubt.  The space will
always be Hausdorff and regular, though.

Thus, for example, the space of continuous loops is a smooth manifold,
as is the space of \(H^1\)-Sobolev loops.  However, this technique is
not going to define the space of \(L^2\)-loops.

\newpage

\section{Looping Bundles}
\label{sec:bundles}

In this section we consider the general construction of looping a
bundle.  That is, given a bundle \(\pi : E \to M\) we consider the
resulting triple \(\pi^L : L E \to L M\).  Our central theme is vector
bundles, but we take in principal bundles, gauge bundles, and
connections on the way.  We start by considering the tangent space of
\(L M\), showing that it is diffeomorphic to the loop of the tangent
space of \(M\).  We conclude by commenting that this does not hold for
the cotangent bundle.

\subsection{The Tangent Space}

Tangent spaces in infinite dimensions can be problematic: there is the
\emph{kinematic} tangent space consisting of the derivatives of short
curves and the \emph{operational} tangent space consisting of
derivations of functions.  In finite dimensions these are the same,
but that need not hold in infinite dimensions.  The kinematic always
satisfies \(T E \cong E \times E\) for a convenient vector space
\(E\), but the operational tangent space may be much larger -- it at
least contains the bidual (\cite[28.3]{akpm}).  We shall not go into a
detailed discussion here since for \(L \R^n\), \cite[Theorem
28.7]{akpm} implies that the two notions coincide.  For more on the
types of tangent vector see \cite[\S 28]{akpm}.

We start with a straightforward result on the structure of the tangent
bundle.

\begin{proposition}
  The tangent bundle of \(L M\), \(T L M\), has the structure of a
  bundle of \(L \R\)--modules.
\end{proposition}

\begin{proof}
  We have local charts \(U_\alpha \cong L \R^n\) and so locally \(T
  U_\alpha \cong L \R^n \times L \R^n\).  Thus we can attempt to
  transfer the natural \(L \R\)--module structure on \(L \R^n\) to the
  tangent spaces of \(L M\).  Lemma~\ref{lem:deriv} shows that we can
  do this as from it we deduce that the derivatives of the transition
  functions are \(L \R\)--linear.
\end{proof}

Our main theorem is the following characterisation of the tangent
space.

\begin{theorem}
  \label{th:tanspace}
  There is a natural diffeomorphism \(T L M \to L T M\) covering the
  identity on \(L M\).
\end{theorem}

The idea is that a small perturbation of a loop involves a small
perturbation of each of its points and so we have a loop of small
perturbations of points.

\begin{proof}
  The evaluation map \(S^1 \times L M \to M\) differentiates to a
  smooth map: \(T S^1 \times T L M \to T M\).  We compose this with
  the the zero section \(S^1 \to T S^1\) to get a smooth map \(S^1
  \times T L M \to T M\).  This is the adjoint of a map \(T L M \to L
  T M\) which is thereby shown to be smooth.  As the map \(S^1 \times
  T L M \to T M\) projects down to the evaluation map \(S^1 \times L M
  \to M\), the map \(T L M \to L T M\) projects down to the identity
  on \(L M\).

  In the case of Euclidean space, we have a natural identifications of
  \(T L \R^n\) with \(L \R^n \times L \R^n\) and of \(T \R^n\) with
  \(\R^n \times \R^n\).  It is simple to show that under these
  identifications, the map \(S^1 \times T L \R^n \to T \R^n\)
  described above becomes \((t, \alpha, \beta) \to (\alpha(t),
  \beta(t))\).  The adjoint of this is the canonical isomorphism \(L
  \R^n \times L \R^n \to L (\R^n \times \R^n)\), which is a
  diffeomorphism.  Hence the result is true for Euclidean space.

  For the more general case we consider the local situation.  We
  choose a chart at a loop \(\alpha \in L M\).  Using a trivialisation
  of \(\alpha^* T M\) we extend the chart map to \(\Psi_\alpha : L
  \R^n \to U_\alpha\).  This trivialisation also extends
  diagram~\eqref{diag:eval} of the proof of theorem~\ref{th:expmfd} in
  the obvious way.  Hence the evaluation map  \(S^1 \times U_\alpha
  \to M\) factors as:
  \[
  S^1 \times U_\alpha \xrightarrow{1 \times {\Psi_\alpha}^{-1}} S^1
  \times L \R^n \to S^1 \times \R^n \to T M \xrightarrow{\eta} M,
  \]
  We differentiate this and evaluate on the zero section to get:
  \[
  S^1 \times T U_\alpha \to S^1 \times T L \R^n \to S^1 \times T \R^n
  \to T^{(2)} M \xrightarrow{d \eta} T M.
  \]
  We now 'de-adjointise' this.  As the second map is \((t, \alpha) \to
  (t, \alpha(t))\), the image of \(T L \R^n\) in \(L(S^1 \times T
  \R^n)\) is \(\Gamma_{S^1}(S^1 \times T \R^n)\) which we identify in
  the usual way with \(L T \R^n\).  Thus we get:
  \[
  T U_\alpha \to T L \R^n \to L T \R^n \to L T^{(2)} M
  \to L T M.
  \]
  From proposition~\ref{prop:charttm}, the induced map \(L T \R^n \to
  L T M\) is a chart map, hence a diffeomorphism onto its image.  The
  first map is the derivative of a chart map for \(L M\), hence is a
  diffeomorphism.  We have already shown that the second map is a
  diffeomorphism.

  Thus the map \(T L M \to L T M\) is a local diffeomorphism where by
  ``local'' we mean in \(L M\).  As both are spaces over \(L M\) and
  the map is of spaces over \(L M\) we therefore have a
  diffeomorphism.
\end{proof}

\begin{corollary}
  \label{cor:tanfibre}
  There is a canonical identification as \(L \R\)--modules of
  \(T_\alpha L M\) with \(\Gamma_{S^1}(\alpha^* T M)\).
\end{corollary}

\begin{proof}
  As the isomorphism \(T L M \to L T M\) covers the identity on \(L
  M\), the tangent space at \(\alpha \in L M\) identifies with the set
  of loops in \(T M\) which project down to \(\alpha\).  As we have
  already seen, this is naturally \(\Gamma_{S^1}(\alpha^* T M)\).

  For the \(L \R\)--module structure, we observe that a trivialisation
  of \(\alpha^* T M\) identifies \(T_\alpha L M\) with \(T_0 L \R^n\)
  and \(\Gamma_{S^1}(\alpha^* T M)\) with \(L \R^n\).  The
  identification between these is the canonical identification of
  \(T_0 L \R^n\) with \(L \R^n\) which is \(L \R\)--linear.
\end{proof}

This needs a word of explanation.  The chart map at \(\alpha \in L
M\), \(\Psi_\alpha : \Gamma_{S^1}(\alpha^* T M) \to U_\alpha\)
identifies the tangent space at a point of \(U_\alpha\) with
\(\Gamma_{S^1}(\alpha^* T M)\) and so in particular we have an
identification of the tangent space at \(\alpha\) with
\(\Gamma_{S^1}(\alpha^* T M)\).  The problem with this is that it may
depend on the local addition \(\eta\) used to define the chart map.
Moreover, we have also relaxed the definition of a local addition so
that the vector bundle did not need to be \(T M\), though \emph{a
fortiori} it has to lie in the same isomorphism class.  In this case,
the chart map identifies the tangent space of \(\alpha\) with
\(\Gamma_{S^1}(\alpha^* E)\) which, although isomorphic to
\(\Gamma_{S^1}(\alpha^* T M)\) does not seem to be canonically so.

The solution is that differentiating the diffeomorphism \((\pi \times
\eta) : E \to V\) defines an isomorphism of bundles \(\pi^* E \cong
\ker d \pi \to (\pi \times \eta)^* T M\), where the \(T M\) is the
second factor in \(T (M \times M)\).  Over the zero section of \(E\)
we thus get an isomorphism \(E \to T M\) defined canonically from the
local addition.  The isomorphism \(T_\alpha L M \to
\Gamma_{S^1}(\alpha^* T M)\) uses this isomorphism of bundles.  One
caveat of this is that if it so happens that our local addition was
defined on \(T M\) then the canonical identification of the above
corollary may not be the simple identification coming from the chart
map.

Now that we have this identification of the tangent bundle, we can
extend proposition~\ref{prop:loopmap}.

\begin{proposition}
  \label{prop:tanloopmap}
  Let \(f : M \to N\) be a smooth map between finite dimensional
  smooth manifolds.  Under the identifications \(T L M \cong L T M\)
  and \(T L N \cong L T N\) we have: \(d (f^L) = (d f)^L\).
\end{proposition}

\begin{proof}
  This is a local result so it is sufficient to work in \(L \R^n\),
  where it is obvious.
\end{proof}

\subsection{Vector Bundles}
\label{sec:vect}

We start by considering various loop spaces associated to a vector
bundle.  We shall formulate our results for a general finite
dimensional vector bundle \(\xi := E \xrightarrow{\pi_\xi} M\).  To
avoid dealing with the twisted situation, we assume that \(E\) is
orientable.

The fibrewise vector space structure on \(E\) can be thought of as
consisting of two maps: \(\R \times E \to E\) and \(E \times_M E \to
E\), satisfying certain relations.  As before, \(E \times_M E\) is an
embedded submanifold of \(E \times E\) and so a map into \(E \times_M
E\) is smooth if and only if the composition with the two projections
to \(E\) are smooth.

We can loop both of these maps and get: \(L \R \times L E \to L E\)
and \(L (E \times_M E) \to L E\).  These satisfy similar relations to
the ones satisfied by the maps on \(E\).  To identify \(L (E \times_M
E)\), from the remark above this consists of pairs of loops in \(E\)
which, time-by-time, project down to the same point on \(M\).  That
is:
\[
L (E \times_M E) = \{(\alpha, \beta) \in L E \times L E :
\pi_\xi(\alpha(t)) = \pi_\xi(\beta(t)) \text{ for all } t \in S^1\}.
\]
Now \(\pi_\xi(\alpha(t)) = \pi_\xi(\beta(t))\) for all \(t\) if and
only if \({\pi_\xi}^L(\alpha) = {\pi_\xi}^L(\beta)\).  Thus \(L (E
\times_M E) = L E \times_{L M} L E\).

Hence we have a structure on the fibres of \(\pi_{L \xi} :=
{\pi_\xi}^L : L E \to L M\) similar to that on the fibres of \(\pi_\xi
: E \to M\).  That is, the fibres of \(L E \to L M\) naturally have
the structure of \(L \R\)-modules.  Note that this implies that they
are vector spaces as there is an inclusion \(\R \to L \R\).  What we
wish to show is that this structure is locally trivial and thus
defines an \(L \R\)-module bundle over \(L M\).

We start by noting that, as for the case \(E = T M\), the fibre of \(L
E\) above a loop \(\alpha \in L M\) can be naturally considered to be
\(\Gamma_{S^1}(\alpha^* E)\).  Moreover, the \(L \R\)--module
structure defined above on the fibres of \(L E \to L M\) is the
natural \(L \R\)--module structure on \(\Gamma_{S^1}(\alpha^* E)\).
Thus once we have established that \(L E\) is a locally trivial bundle
of \(L \R\)--modules over \(L M\), the isomorphism of
theorem~\ref{th:tanspace} will be of \(L R\)--module bundles.

The mainstay of the analysis of the loop space of \(E\) in terms of
that of \(M\) is that we do not need to use a local addition on the
whole of \(E\) to define the charts for \(L E\).  All we need is the
local addition on \(M\) together with a connection on \(E\).

The bundle \(T M \oplus E\) can be realised by the \emph{Whitney sum}.
This identifies the total space of \(T M \oplus E\) with the pull-back
of \(T M \times E \to M \times M\) by the diagonal map \(M \to M
\times M\).  Thus the total space of \(T M \oplus E\) is:
\[
T M \times_M E := \{(u,v) \in T M \times E : \pi(u) = \pi_\xi(v)\}.
\]
The obvious projection maps \(T M \times_M E\) to \(T M\) and to \(E\)
fit together to give a commuting diagram:
\[
\begin{CD}
T M \times_M E @>>> E \\
@VVV @VV\pi_\xi V \\
T M @>\pi>> M.
\end{CD}
\]
Thus we can think of \(T M \times_M E\) as the total space of any of:
\(T M \oplus E \to M\), \(\pi^* E \to T M\), or \(\pi_\xi^* T M \to
E\).

Let \(\eta : T M \to M\) be a local addition on \(M\) with \(V
\subseteq M \times M\) the corresponding open neighbourhood of the
diagonal.  For \(v \in T M\), define \(P_u : E_{\pi(u)} \to
E_{\eta(u)}\) to be the operation of parallel transport along the path
\(t \to \eta(t u)\) (recall that \(\eta(0 u) = \pi(u)\)).  Using this,
we define a variant of local addition.  Define \(\mu : T M \times_M E
\to E\) by \(\mu(u,v) = P_u(v)\).  Note that this is well-defined
since \(\pi_\xi(v) = \pi(u)\) so \(v \in E_{\pi(u)}\).  Observe that
\(\pi_\xi \mu (u,v) = \eta(u)\) so we have a commutative diagram:
\[
\begin{CD}
T M \times_M E @>\mu>> E \\
@VVV @VV\pi_\xi V \\
T M @>\eta >> M.
\end{CD}
\]

Consider the map \(\pi \times \mu : T M \times_M E \to M \times E\)
where, by abuse of notation, the map \(\pi : T M \times_M E \to M\) is
\((u,v) \to \pi(u) = \pi_\xi(v)\).  Let \(E_V := \{(x,v) \in M \times
E : (x, \pi_\xi(v)) \in V\}\).

We claim that \(\pi \times \mu : T M \times_M E \to E_V\) is a
diffeomorphism.  The most geometric way to see this is to observe that
there is another space hidden in the background, namely:
\[
T M \times_\eta E := \{(u,v) \in T M \times E : \eta(u) =
\pi_\xi(v)\}.
\]
This is the pull-back of \(E\) over \(T M\) via \(\eta\), so the
following diagram commutes:
\[
\begin{CD}
T M \times_\eta E @>\tilde{\mu}>> E \\
@VVV @VV\pi_\xi V \\
T M @>\eta>> M
\end{CD}
\]
The map \(\pi \times \eta : T M \to V \subseteq M \times M\) is a
diffeomorphism, so \(\pi \times \tilde{\mu} : T M \times_\eta E \to
E_V \subseteq M \times E\) is also a diffeomorphism.

Now \(T M \times_\eta E\) is isomorphic to \(T M \times_M E\) as
bundles over \(T M\) since \(\pi : T M \to M\) and \(\eta : T M \to
M\) are homotopic maps.  An explicit homotopy is \(h(v,t) = \eta(t
v)\).  This defines an explicit isomorphism between the bundles by
parallel transporting the bundle along the fibres of the homotopy,
i.e.~along the paths \(t \to \eta(t v)\).  This isomorphism is \((u,v)
\to (u, P_u(v))\).  Hence the composition of this with \(\tilde{\mu}\)
is \(\mu\).  Thus \(\pi \times \mu : T M \times_M E \to E_V\) is a
diffeomorphism.

We wish to play the same game for \(\mu : T M \times_M E \to E\) that
we had for \(\eta : T M \to M\).  That is, we wish to construct maps
that have the effect of sending a section \(\beta\) of something to
the loop \(t \to \mu(\beta(t))\) in \(E\).  Clearly, we must be able
to think of \(\beta\) as a loop in \(T M \times_M E\), just as in the
case of \(M\) we thought of a section of \(\alpha^* T M\) as a loop in
\(T M\).  We have three candidates for what the section is of since
\(T M \times_M E\) has the structure of a vector bundle over each of
\(T M\), \(E\), and \(M\).  To determine the correct one, we recall
from the story for \(M\) that the two components of the map \(\pi
\times \eta : T M \to V \subseteq M \times M\) were used as follows:
the second defined the chart map, the first told us which chart we
were in.

Therefore, as we have \(\pi \times \mu : T M \to E_V \subseteq M
\times E\), we should try to index the charts by loops in \(M\),
\emph{not} by loops in \(E\).  Thus the domain of the chart map will
be sections of \(\alpha^*(T M \times_M E)\) for \(\alpha \in L M\).
We can write this as \(\alpha^* (T M \oplus E)\) to indicate this
without constantly needing to remind ourselves that \(\alpha\) is a
loop in \(M\).

The codomain of the chart map is the set of loops in \(E\) with the
property that \((\alpha(t), \beta(t)) \in E_V\), for all \(t \in
S^1\).  Now \(\beta\) satisfies this condition if and only if
\((\alpha(t), \pi_\xi(\beta(t))) \in V\) for all \(t \in S^1\).  Thus
the codomain is \(\{\beta \in L E : {\pi_\xi}^L(\beta) \in
U_\alpha\}\).  Write this as \(L E_{U_\alpha}\).  Note that this is
the pull-back of \(L E \to L M\) under the inclusion \(U_\alpha \to L
M\).  These cover \(L E\) because the \(U_\alpha\) cover \(L M\).

One can alter \(\mu : T M \times_M E \to E\) to a full local addition
on \(E\) by taking the pull-back of \(T M \oplus E \to M\) via
\(\pi_\xi : E \to M\).  The resulting bundle is isomorphic to \(T E\)
and the map \(\mu {\pi_\xi}_* : \pi_\xi^*(T M \oplus E) \to E\) is a
local addition.  This can then be used to show that the smooth
structure coming from the charts constructed above is the same as the
natural smooth structure on \(L E\).

We return to the charts defined by \(\mu : T M \times_M E \to E\).
The domain of a chart is \(\Gamma_{S^1}(\alpha^* (T M \oplus E))\).
The total space of \(\alpha^*(T M \oplus E)\) is the subset of \(S^1
\times (T M \times_M E)\) such that \(\alpha(t) = \pi(u,v)\).  From
the definition of \(T M \times_M E\) we can see that this is the
subset of \(S^1 \times T M \times E\) such that \(\alpha(t) = \pi(u) =
\pi_\xi(v)\).  Therefore a section of \(\alpha^* (T M \oplus E)\) can
be written in the form \(t \to (t, \beta(t), \gamma(t))\) where
\(\beta : S^1 \to T M\) and \(\gamma : S^1 \to E\) are such that
\(\alpha = \pi^L\beta = {\pi_\xi}^L\gamma\) for all \(t \in
S^1\).  Thus the map:
\[
\big(t \to (t, \beta(t), \gamma(t)) \big) \to \Big( \big( t \to (t,
\beta(t)) \big), \big( t \to (t, \gamma(t)) \big) \Big)
\]
identifies \(\Gamma_{S^1}(\alpha^*(T M \oplus E))\) with
\(\Gamma_{S^1}(\alpha^* T M) \times \Gamma_{S^1}(\alpha^* E)\).
(It is simple to show that the identification \(L \R^{n + m} \cong
L \R^n \times L \R^m\) preserves all the structure including the
smooth structure, so the identification of the spaces of sections
given above also preserves the smooth structure.)

We have:
\[
\Theta_\alpha : \Gamma_{S^1}(\alpha^* T M) \times
\Gamma_{S^1}(\alpha^* E) \to L E_{U_\alpha}
\]
satisfying:
\[
\Theta_\alpha(\beta, \gamma)(t) = \mu(\tilde{\beta}(t),
\tilde{\gamma}(t)) =
P_{\tilde{\beta}(t)}(\tilde{\gamma}(t)),
\]
where \(\beta(t) = (t, \tilde{\beta}(t))\) and \(\gamma(t) = (t,
\tilde{\gamma}(t))\).

Let \(\gamma_1, \gamma_2 \in \Gamma_{S^1}(\alpha^* E)\).  Let
\(\tilde{\gamma}_1, \tilde{\gamma}_2 \in L E\) be the
corresponding loops.  Let \(\nu \in L \R\).  The loop in \(L E\)
corresponding to the section \(\gamma_1 + \nu \gamma_2\) is
\(\tilde{\gamma}_1 + \nu \tilde{\gamma}_2\).  Therefore, since
\(v \to P_u(v)\) is linear:
\begin{align*}
  \Theta_\alpha(\beta, \gamma_1 + \nu \gamma_2)(t) &=
  P_{\tilde{\beta}(t)}(\tilde{\gamma}_1(t) + \nu(t)
  \tilde{\gamma}_2(t)) \\
  & = P_{\tilde{\beta}(t)}(\tilde{\gamma}_1(t)) + \nu(t)
  P_{\tilde{\beta}}(\tilde{\gamma}_2(t)) \\
  & = \big(\Theta_\alpha(\beta, \gamma_1) + \nu
  \Theta_\alpha(\beta,\gamma_2) \big)(t).
\end{align*}
Hence \(\Theta_\alpha\) is \(L\R\)-linear in its second argument.

Thus \(L E \to L M\) is a locally trivial \(L\R\)-module bundle.  The
fibres are modelled on the spaces \(\Gamma_{S^1}(\alpha^* E)\).  Since
we assumed that \(E\) is orientable, these are all isomorphic to \(L
\R^n\) (as \(L \R\)-modules) via a smooth trivialisation of \(\alpha^*
E \to S^1\).

\subsection{Principal and Gauge Bundles}

The key in the analysis of the vector bundle was the existence of the
parallel transport operator and the fact that it is linear.  There are
other bundles with such operators, the most common being principal
bundles.

Let \(G\) be a finite dimensional Lie group, \(Q \to M\) a principal
\(G\)-bundle.  As always, we shall assume that it is \emph{orientable}
in that \(\alpha^* Q \to S^1\) admits a section for any loop \(\alpha
: S^1 \to M\).  Now \(L Q\) inherits an action of the loop group \(L
G\) by \((\alpha \cdot \gamma)(t) = \alpha(t) \cdot \gamma(t)\).  The
space of sections of \(\alpha^*Q\) also has such an action.  A
trivialisation, \(\alpha^* Q \cong S^1 \times G\), identifies
\(\Gamma_{S^1}(\alpha^* Q)\) with \(L G\) and this is an isomorphism
of \(L G\)-spaces.

A connection on \(Q\) defines a parallel transport operator along any
path and this preserves the fibrewise \(G\)-action.  Therefore the
same analysis as for the case of the vector bundle leads to maps:
\[
\Theta_\alpha : \Gamma_{S^1}(\alpha^* T M) \times
\Gamma_{S^1}(\alpha^* Q) \to L Q_{U_\alpha}
\]
which are \(L G\)-equivariant.  This shows that \(L Q \to L M\) is a
locally trivial principal \(L G\)-bundle.

Suppose that \(S\) is a (finite dimensional) smooth manifold with a
\(G\)-action.  We can form the locally trivial bundle \(R := Q
\times_G S \to M\).  It is clear that \(L R\) can be described as \(L
Q \times_{L G} L S \to L M\) and that this is a locally trivial bundle
modelled on \(L S\).  Structure on \(S\) that is preserved by \(G\)
defines structure on \(R\) so since \(L G\) preserves the ``looped
structure'' on \(L S\), there is a corresponding structure on \(L R\).
We have already seen this in the vector bundle situation in that the
\(\R\)-module (i.e.~vector space) structure of the fibres of \(E\)
looped to give an \(L \R\)-module structure on the fibres of \(L E\).

Another important situation where this arises is when \(G\) acts on
itself via the adjoint action.  This defines the \emph{gauge group}
corresponding to \(Q\), \(Q^\ad := Q \times_\ad G\) which is a bundle
of groups over \(M\) modelled on \(G\).  The corresponding loop space,
\(L Q^\ad\), is a bundle of groups over \(L M\) modelled on \(L G\).
We have the identity \(L (Q^\ad) = (L Q)^\ad := L Q \times_{\ad} L
G\).

Returning to the vector bundle case, if \(E = Q \times_G \R^n\) then
\(L E = L Q \times_{L G} L \R^n\), as \(L \R\)-modules.  Now as \(Q\)
is the frame bundle of \(E\), it would be nice to be able to say that
\(L Q\) is the frame bundle of \(L E\).  This is not true as a
statement about vector bundles, since \(\gl(L \R^n) \ne L \gl(\R^n)\),
but is true as a statement about \(L\R\)-modules as a consequence of:

\begin{proposition}
  Let \(g : L \R^n \to L \R^n\) be an isomorphism of \(L
  \R\)-modules.  Then \(g \in L \gl(\R^n)\).
\end{proposition}

\begin{proof}
  For a ring \(R\), module \(M\), and subset \(S \subseteq M\), let
  \(\langle S \rangle_R\) denote the \(R\)-linear span of \(S\) in
  \(M\).  We add the subscript \(R\) to the usual notation as we will
  have a space being a module over two different rings.

  It is mildly obvious that the space \(L \R^n\) is a free \(L
  \R\)-module of rank \(n\).  Let \(\m{B}\) be an \(n\)-element
  generating set (basis).

  Let \(t \in S^1\).  The evaluation map \(e_t : L \R \to \R\) is a
  ring map and so converts any \R-module (i.e.~vector space) into an
  \(L \R\)-module.  The corresponding evaluation map \(e_t : L \R^n
  \to \R^n\) is then a map of \(L \R\)-modules.  Then:
  \[
  \R^n = e_t(L \R^n) = \langle e_t(\m{B}) \rangle_{L \R} = \langle
  e_t(\m{B}) \rangle_\R.
  \]
  The first equality is because \(e_t : L \R^n \to \R^n\) is
  surjective and the last because the action of \(L \R\) on \(\R^n\)
  factors through \R.  As \(e_t(\m{B})\) is an \(n\)-element set, it
  is therefore a basis for \(\R^n\).

  Let \(g : L \R^n \to L \R^n\) be an \(L \R\)-isomorphism.  For \(j =
  1, \dotsc, n\) let \(\alpha_j \in L \R^n\) be the constant loop at
  the \(j\)th element of the standard basis of \(\R^n\).  Define a
  loop \(\gamma\) in the space of \(n \times n\) matrices,
  \(M_n(\R)\), by arranging the loops \(g \alpha_j\) in columns:
  \[
  \gamma(t) = \big( (g \alpha_1)(t), \dotsc, (g \alpha_n)(t) \big).
  \]

  Now \(\{\alpha_1, \dotsc, \alpha_n\}\) is a generating set for \(L
  \R^n\) as an \(L \R\)-module.  As \(g\) is an \(L \R\)-isomorphism,
  \(\{g \alpha_1, \dotsc, g \alpha_n\}\) is also a generating set for
  \(L \R^n\).  Therefore at each \(t \in S^1\), \(\{(g \alpha_1)(t),
  \dotsc, (g \alpha_n)(t)\}\) is a basis for \(\R^n\).  Hence
  \(\gamma\) takes values in \(\gl(\R^n)\).  

  Now the action of \(L \gl(\R^n)\) on \(L \R^n\) is \(L \R\)-linear
  as it is the loop of the \R-linear action of \(\gl(\R^n)\) on
  \(\R^n\).  As \(\alpha_j\) is the constant loop at the \(j\)th
  standard basis element, \(\gamma \alpha_j\) selects the \(j\)th
  column of \(\gamma\), which is \(g \alpha_j\).  Hence \(g\) and
  \(\gamma\) are both \(L \R^n\)-linear maps which agree on a
  generating set, thus are the same map.
\end{proof}

The proof of this result shows that the correct interpretation of an
\(L \R\)-frame of \(L E \to L M\) at a loop \(\alpha \in L M\) is a
smooth choice of frame of each \(E_{\alpha(t)}\).  Thus an \(L
\R\)-frame of \(L E\) is a trivialisation of \(\alpha^* E\).  To get
back the isomorphism \(L \R^n \to L_\alpha E\), we use the
trivialisation to identify \(L \R^n\) with \(\Gamma_{S^1}(\alpha^*
E)\) which is naturally identified with the fibre \(L_\alpha E\).

\subsection{Connections}
\label{sec:connect}

In this section we consider connections on vector bundles.  As with
the vector bundles themselves, we are interested in structure on the
loop space that arises from structure on the original space.  Our main
theorem in this section is:

\begin{theorem}
  \label{th:connection}
  A principal connection in finite dimensions loops to a principal
  connection on the loop space.
\end{theorem}

There are several ways of thinking of connections in finite dimensions
and these carry over to infinite dimensions.  The general theory is
contained in \cite[\S 37]{akpm} of which we recall the basics here.

Let \(G\) be a Lie group, \(\pi : P \to M\) a principal \(G\)--bundle.
We are considering both the finite and infinite dimensional cases
here.  The two main ways that we think of a connection are:
\begin{enumerate}
\item A projection \(\Phi : T P \to V P\), where \(V P = \ker d \pi\)
  is the vertical tangent bundle, which is \(G\)--equivariant for the
  \(G\) action on \(P\).  That is, for each \(g \in G\) and \(u \in
  P\) the following diagram commutes:
  \[
  \begin{CD}
    T_u P @>d r^g>> T_{u g} P \\
    @V\Phi VV @V\Phi VV \\
    V_u P @>d r^g>> V_{u g} P,
  \end{CD}
  \]
  where \(r^g : P \to P\), \(r^g(u) = u g\), is the action of \(g \in
  G\).

\item The connection one-form \(\omega : T P \to \mf{g}\), the Lie
  algebra of \(G\), satisfying:
  \begin{enumerate}
  \item \(w(\zeta_X(u)) = X\) for all \(X \in \mf{g}\), where
    \(\zeta_X \in \mf{X}(P)\) is the vector field defined by
    \(\zeta_X(u) := d r_{(u,e)}(0, X)\) (note that \(d r_{(u,e)}\) is
    a linear map \(T_u P \times T_e G \to T_u P\));

  \item \(w\) is \(G\)--equivariant: \(((r^g)^*\omega)(X) =
    \Ad_{g^{-1}} \omega(X)\) for all \(g \in G\) and \(X \in T_u P\).
  \end{enumerate}
\end{enumerate}
The equivalence is given by \(\Phi(X) = \zeta_{\omega(X)}(u)\) for \(X
\in T_u P\).

We can loop both maps to get \(\Phi^L : L T P \to L V P\) and
\(\omega^L : L T P \to L \mf{g}\).  Using theorem~\ref{th:tanspace} we
can identify \(L T P\) with \(T L P\) in both cases.  To proceed
further, we need some technical results:

\begin{lemma}
  Let \(G\) be a finite dimensional Lie group, \(\pi : P \to M\) a
  principal \(G\)--bundle over a finite dimensional manifold.  Then:
  \begin{enumerate}
  \item The Lie algebra of \(L G\) is canonically isomorphic to \(L
    \mf{g}\).  Although we shall not need it here, it seems an
    appropriate place to mention that the exponential map \(L \mf{g}
    \to L G\) is the loop of the exponential map \(\mf{g} \to G\).

  \item Under the isomorphism \(T L P \cong L T P\), the vertical
    tangent bundle \(V L P\) is mapped to \(L V P\).

  \item Under the isomorphism \(T L P \cong L T P\), the vector field
    mapping \(\zeta : L \mf{g} \to \mf{X}(L P)\) is the loop of
    \(\zeta : \mf{g} \to \mf{X}(P)\) followed by the canonical
    inclusion \(L \mf{X}(P) \to \Gamma_{L P}(L T P) \cong \mf{X}(L
    P)\) which is given by:
    \(
    \chi_\alpha(t) = \chi(t)_{\alpha(t)}.
    \)

  \item The adjoint action of \(L G\) on \(L \mf{g}\) is the loop of
    the adjoint action of \(G\) on \(\mf{g}\).

  \end{enumerate}
\end{lemma}

\begin{proof}
  \begin{enumerate}
  \item The Lie algebra of any Lie group is the (kinematic) tangent
    space at the identity.  Since the identity of \(L G\) is the
    constant loop at the identity, \(e\), of \(G\), the Lie algebra of
    \(L G\) is, by corollary~\ref{cor:tanfibre}, naturally identified
    with \(\Gamma_{S^1}({\gamma_e}^* T G)\).  As \(\gamma_e\) is
    constant at \(e\), \({\gamma_e}^* T G = S^1 \times T_e G\) and so
    \(T_{\gamma_e} L G\) is naturally identified with loops in \(T_e
    G\), i.e.~with \(L \mf{g}\).

    The statement about the exponential maps follows because the
    evaluation maps, \(e_t : L G \to G\), are group homomorphisms.
    The induced Lie algebra map is also the evaluation map, \(e_t : L
    \mf{g} \to \mf{g}\).  Thus for \(t \in S^1\), \(s \in \R\), and
    \(\chi \in L \mf{g}\), \(\exp(s \chi)(t) = \exp(s \chi(t))\).

  \item We need to identify \(\ker d (\pi^L) : T L P \to T L M\) under
    the equivalence \(T L P \cong L T P\) and \(T L M \cong L T M\).
    From proposition~\ref{prop:tanloopmap}, \(d (\pi^L)\) is \((d
    \pi)^L\).  Thus \((d \pi)^L \beta = 0\) if and only if \(d \pi
    \beta(t) = 0\) for all \(t\) and so \(\ker (d \pi)^L \subseteq L T
    P\) is the loop of \(\ker d \pi \subseteq T P\), i.e.~\(V L P\)
    corresponds to \(L V P\).

  \item Let \(\chi \in L \mf{g}\) and \(\alpha \in L P\).  By
    definition:
    \[
    \zeta_\chi(\alpha) = d r_{(\alpha, \gamma_e)}(0, \chi) :
    T_{\alpha} L P \times T_{\gamma_e} L G \to T_{\alpha} L P.
    \]
    As the action of \(L G\) on \(L P\) is the loop of the action of
    \(G\) on \(P\), the various isomorphisms mean that: \(d
    r_{(\alpha, \gamma_e)}(0, \chi)(t) = d r_{(\alpha(t), e)}(0,
    \chi(t))\).  Hence:
    \[
    \zeta_\chi(\alpha)(t) = \zeta_{\chi(t)}(\alpha(t))
    \]
    as required.

  \item The adjoint action of \(L G\) on \(L \mf{g}\) is the
    derivative of the conjugation action at the identity.  The
    conjugation action is the loop of the conjugation action of \(G\)
    on \(\mf{g}\).  Hence its derivative is the loop of the
    derivative. \qedhere
  \end{enumerate}
\end{proof}

We can now prove theorem~\ref{th:connection}.

\begin{proof}[Proof of theorem~\nmref{th:connection}]
  Let \(G\) be a finite dimensional Lie group and let \(P \to M\) be a
  principal \(G\)--bundle over a finite dimensional manifold.  Let
  \(\Phi : T P \to V P\) be a principal connection on \(P\).  Using
  the equivalences \(T L P \cong L T P\) and \(V L P \cong L V P\), we
  loop \(\Phi\) to a map \(\Phi^L : T L P \to V L P\) which we claim is
  a principal connection on \(L P\).  It is certainly a fibre
  projection as at a loop \(\alpha \in L P\) it is the projection
  \(\Gamma_{S^1}(\alpha^* T P) \to \Gamma_{S^1}(\alpha^* V P)\).  It
  is also \(L G\)--equivariant as \(\Phi : T P \to V P\) is
  \(G\)--equivariant and the action is defined pointwise.

  For the other view of a connection, let \(\omega : T P \to \mf{g}\)
  be the connection form associated to \(\Phi\).  Looping this, using
  the equivalences, yields \(\omega^L : T L P \to L\mf{g}\).  For \(\chi
  \in L \mf{g}\) and \(\alpha \in L P\) we have, for all \(t \in S^1\):
  \[
  \omega^L(\zeta_\chi(\alpha))(t) = \omega(\zeta_\chi(\alpha)(t)) =
  \omega(\zeta_{\chi(t)}(\alpha(t))) = \chi(t)
  \]
  and hence \(\omega^L(\zeta_\chi(\alpha)) = \chi\).  Secondly,
  \(\omega^L\) is \(L G\)--equivariant as all the actions are loops of
  the actions on the original spaces.

  Finally, the equation \(\Phi(X) = \zeta_{\omega(X)}(u)\) for \(X \in
  T_u P\) shows that for \(\chi \in T_\alpha L P\):
  \[
  \Phi(\chi)(t) = \Phi(\chi(t)) = \zeta_{\omega(\chi(t))}(\alpha(t))
  \]
  and hence \(\Phi^L(\chi) = \zeta_{\omega^L(\chi)}(\alpha)\).
\end{proof}

Associated to a connection is a parallel transport operation that
lifts curves in the base to curves in the fibre such that the
derivative is horizontal (i.e.~in \(\ker \Phi\)).  The parallel
transport operations on \(M\) and on \(L M\) are related as follows:

\begin{proposition}
  \label{prop:parallel}
  The parallel transport operator exists in \(L P\) and corresponds to
  the parallel transport operator in \(P\) under the evaluation maps,
  \(e_t : L P \to P\).  That is, if \(\tilde{c} : \R \to L P\) is a
  horizontal lift of \(c : \R \to L M\) then \(e_t \tilde{c} : \R \to
  P\) is a horizontal lift of \(e_t c : \R \to M\).
\end{proposition}

\begin{proof}
  As a loop is completely determined by the values it takes, this
  characterisation completely specifies the parallel transport
  operator in \(L P\).  As parallel transports are unique providing
  they exist, \cite[37.6]{akpm}, we just need to show that this
  characterisation \emph{defines} the parallel transport.

  Let \(c : \R \to L M\) be a curve and \(v \in L P\) a lift of
  \(c(0)\).  For \(t \in S^1\) let \(c_t := e_t c : \R \to M\) and
  \(v_t := e_t v \in P\), so that \(v_t\) is a lift of \(c_t(0)\).
  Let \(\tilde{c}_t : \R \to P\) be the parallel transport of \(v_t\)
  along \(c_t\).  We therefore have a map \(S^1 \times \R \to P\),
  \((t,s) \to \tilde{c}_t(s)\).  Let \(\tilde{c} : \R \to L P\) be the
  adjoint of this.

  As parallel transport in \(P\) is smooth in all initial conditions,
  the map \((t,s) \to \tilde{c}_t(s)\) is smooth.  Hence \(\tilde{c}\)
  is smooth.  Clearly, \(\pi^L \tilde{c} = c\) and \(\tilde{c}(0) =
  v\).  Finally, \(e_t \Phi^L \tilde{c}' = \Phi {\tilde{c}_t}' = 0\)
  so \(\Phi^L \tilde{c}' = 0\) and hence \(\tilde{c}\) is horizontal.
  Thus it is the parallel transport of \(v\) along \(c\).
\end{proof}

Given an action of \(G\) on a (finite dimensional) manifold \(S\) we
get an associated fibre bundle \(R := P \times_G S \to M\).  The
connection on \(P\) induces a connection on \(R\).  Similarly, the
connection on \(L P\) induces a connection on \(L R\).  That these
correspond is straightforward to deduce.

In particular, for a representation of \(G\) on \(\R^n\) and \(E := P
\times_G \R^n\) we get linear connections on \(E\) and \(L E\).  From
these we construct covariant differentiation operators.  We start with
the canonical isomorphism \(\text{vl}_E : E \times_M E \to V E\),
where \(\text{vl}_E(u_x, v_x)\) is represented by the short curve \(s
\to u_x + s v_x\).  Using this, we define the \emph{connector} \(K_E\)
of the connection \(\Phi_E\) by:
\[
K_E := \text{pr}_2 \circ (\text{vl}_E)^{-1} \circ \Phi_E : T E \to V E
\to E \times_M E \to E.
\]
The covariant derivative on vector bundles is defined, following
\cite[37.28]{akpm}, as follows: for any manifold \(N\), smooth
mapping \(s : N \to E\), and (kinematic) vector field \(X \in
\mf{X}(N)\) we define the \emph{covariant derivative} of \(s\) along
\(X\) by:
\[
\nabla^E_X s := K_E \circ d s \circ X : N \to T N \to T E \to E.
\]

If \(s\) is a lift of some fixed map \(f : N \to M\) then \(\nabla_X^E
s\) is also a lift of the same map and so we get an induced operation
on sections of \(f^* E\).  In particular, taking \(f : M \to M\) to be
the identity we get the usual operator:
\[
\nabla^E : \mf{X}(M) \times \Gamma_M(E) \to \Gamma_M(E).
\]
The reason for taking the more general approach is that it is the
better setting for expressing the relationship between \(\nabla^E\)
and \(\nabla^{L E}\), which is defined analogously, because by
theorem~\ref{th:expmfd} there is an equivalence \(\Ci(N, L E) \cong
\Ci(S^1 \times N, E)\).  Therefore given a map \(s : N \to L E\) we
take its adjoint to get a map \(s^\lor : S^1 \times N \to E\).  A
vector field \(X\) on \(N\) defines one on \(S^1 \times N\) in the
obvious way which we also denote by \(X\).

\begin{theorem}
  \label{th:covadj}
  \(\displaystyle \big(\nabla^{L E}_X s)^\lor = \nabla^E_X (s^\lor)\).
\end{theorem}

\begin{proof}
It is obvious that \(K_{L E}\) is the loop of \(K_E\) since each map
in the definition of \(K_{L E}\) is the loop of the corresponding map
in the definition of \(K_E\).  It is straightforward to show that the
identification \(f \to f^\lor\) of theorem~\ref{th:expmfd} together
with the identification of \(T L M\) with \(L T M\) means that:
\[
(d f)^\lor = d (f^\lor) \circ (\zeta \times 1) : S^1 \times T N \to T
S^1 \times T N \to T M,
\]
where \(\zeta : S^1 \to T S^1\) is the zero section.  Now the vector
field on \(S^1 \times N\) corresponding to \(X\) can be thought of as
\((\zeta \times 1)(1 \times X)\).  Thus
\(\nabla^{E}_X (s^\lor)\) is:
\[
S^1 \times N \xrightarrow{1 \times X} S^1 \times T N
\xrightarrow{\zeta \times 1} T S^1 \times T N \xrightarrow{d (s^\lor)}
T E \xrightarrow{K_E} E.
\]
the central two terms contract to \((d s)^\lor\) so we have:
\[
S^1 \times N \xrightarrow{1 \times X} S^1 \times T N \xrightarrow{(d
  s)^\lor} T E \xrightarrow{K_E} E.
\]
This is the adjoint of:
\[
N \xrightarrow{X} T N \xrightarrow{d s} L T E \xrightarrow{K_E} L E,
\]
which is \(\nabla^{L E}_X s\) as required.
\end{proof}

In particular, the covariant differential operator:
\[
\nabla^{L E} : \mf{X}(L M) \times \Gamma_{L M}(L E) \to \Gamma_{L M}(L
E)
\]
is adjoint to the covariant differential operator on \(e^* E \to S^1
\times L M\) where \(e : S^1 \times L M \to M\) is the evaluation map.

\begin{lemma}
  Fix \(f : N \to L M\) and consider only those maps that are
  lifts of \(f\).  For such sections, the connection \(\nabla^{L E}\)
  is \(L \R\)--linear.
\end{lemma}

\begin{proof}
  We need to restrict our attention to lifts of a fixed map in order
  to have an addition on the maps \(s : N \to L E\) since we need to
  know that \(s_1(x)\) and \(s_2(x)\) end up in the same fibre of \(L
  E \to L M\).

  Standard properties of connections imply \R--linearity so all that
  we need to prove is the correct behaviour under multiplication by an
  element of \(L \R\).  This follows from the pleasant properties of
  the differential of a map under the adjoint mapping.  Namely, that
  for \(s : N \to L E\) with adjoint \(s^\lor : S^1 \times N \to E\),
  the adjoint of \(d s\) is the \(N\)--directional derivative of
  \(s^\lor\).  That is,
  \[
  (d s)^\lor = d (s^\lor) \circ (\zeta \times 1) : S^1 \times T N \to
  T S^1 \times T N \to T E.
  \]
  Now for \(\beta \in L \R\), \((\beta s)^\lor\) is the map \((t,x)
  \to \beta(t) s^\lor(t,x)\).  When taking the \(N\)--directional
  derivative, only the \(s^\lor\)--factor contributes and we find that
  \(d(\beta s)^\lor = \beta d s^\lor\).  Hence, undoing the adjoints,
  \(d (\beta s) = \beta d s\).

  Since \(K_{L E}\) is the loop of \(K_E\), it is \(L \R\)--linear and
  hence \(K_{L E} (\beta d s) = \beta K_{L E} d s\).  Thus
  \(\nabla_{X}^{L E} s\) is \(L \R\)--linear in \(s\).
\end{proof}

Even when \(N\) is a loop space, such as \(L M\) itself, the
connection will not be \(L \R\)--linear in the other variable.  This
is because for \(\gamma \in L \R^n\) and \(\beta \in L \R\), the
vectors \(\gamma\) and \(\beta \gamma\) represent completely different
directions.  It is tempting to view \(\beta \gamma\) as an elaborate
stretch of \(\gamma\) but really they are unrelated vectors.

A connection on \(T M\) loops to one on \(L T M\), whence
to one on \(T L M\).  Thus we can consider its torsion:
\[
\tau(X,Y) := \nabla_X Y - \nabla_Y X - [X, Y].
\]
As in finite dimensions this is a tensor and so is given by a
fibrewise bilinear, skew-symmetric map \(T L M \times T L M \to T L
M\).

\begin{proposition}
  The torsion on \(L M\) is the loop of the torsion on \(M\).
\end{proposition}

\begin{proof}
  We start with a couple of facts about \(L \R^n\):
  \begin{enumerate}
  \item Any loop in \(\R^n\) is the sum of two never-zero loops.

    Let \(\beta : S^1 \to \R^n\).  As \(S^1\) is compact the image of
    \(\beta\) is bounded and so there is some non-zero \(v \in \R^n\)
    such that \(\beta(t) \ne v\) for all \(t\).  Thus \(\beta -
    \gamma_v\) is never zero, where \(\gamma_v\) is the constant loop
    at \(v\).  Hence \(\beta\) is the sum \((\beta - \gamma_v) +
    \gamma_v\) of two never-zero loops.

  \item If \(n \ge 2\), given a pair of never-zero loops \(\beta\) and
    \(\delta\), there are never-zero loops \(\beta_1\), \(\beta_2\),
    \(\delta_1\), \(\delta_2\) such that: \(\beta = \beta_1 +
    \beta_2\) and \(\delta = \delta_1 + \delta_2\) and, for all \(t\)
    and for \(i,j \in \{1,2\}\), \(\{\beta_i(t), \delta_j(t)\}\) is a
    linearly independent set.

    If \(n \ge 3\) this is simple: choose a never-zero loop \(\gamma\)
    such that \(\gamma(t)\) is not in the linear span of \(\{\beta(t),
    \delta(t)\}\).  Let \(\beta_1 = \frac12(\beta + \gamma)\),
    \(\beta_2 = \frac12(\beta - \gamma)\), \(\delta_1 = \delta_2 =
    \frac12 \delta\).  Then \(\beta_i\) and \(\delta_i\) satisfy the
    required properties.

    If \(n = 2\) this is slightly more complicated.  The idea relies
    on the fact that if the images of two loops lie on the same side
    of a line through the origin then at no time can they be
    colinear.  Thus we arrange matters so that the four loops lie in
    the four quadrants with the \(\beta_i\) in opposite quadrants
    (hence also the \(\delta_i\) in opposite quadrants).  Then each
    pair \(\{\beta_i, \delta_j\}\) lie on the same side of one of the
    axes and so can never be colinear.  Specifically, let \(R > 0\) be
    such that the images of \(\beta\) and \(\delta\) lie within the
    open square with vertices at \((\pm R, \pm R)\).  Let \(\gamma_1\)
    and \(\gamma_2\) be constant loops at two adjacent vertices of
    this square.  Let \(\beta_1 = \beta - \gamma_1\), \(\beta_2 =
    \gamma_1\), \(\delta_1 = \delta - \gamma_2\), \(\delta_2 =
    \gamma_2\).  Then these satisfy the requirements.

  \item Hence if \(\rho_i : L \R^n \times L \R^n \to L \R^n\), \(i =
    1, 2\), are bilinear and such that \(\rho_1(\beta, \gamma) =
    \rho_2(\beta, \gamma)\) whenever \(\beta, \gamma\) are never-zero
    and never-colinear then \(\rho_1 = \rho_2\).
  \end{enumerate}

  These extend to sections of an orientable vector bundle over \(S^1\)
  via a trivialisation since the concepts of ``never-zero'' and
  ``never-colinear'' are preserved under bundle maps.  Hence they
  transfer to the fibres of \(T L M \to L M\).  Note that we really
  need orientable here rather than it being a convenient
  simplification as there are no never-zero sections of the M\"obius
  line.

  The other information that we need to know about is how the
  covariant derivative on \(L E\) -- on any vector bundle --
  transforms under maps of the source.  As in finite dimensions, the
  covariant derivative of \(s : N \to L E\) along \(X \in \mf{X}(N)\)
  at \(x \in N\) depends on \(X\) only up to \(X(x)\).  That is, the
  map \(\nabla_{X(x)} : \Ci(N, L E) \to L E\) makes sense and
  satisfies:
  \begin{align*}
    \pi(\nabla_{X(x)} s) &= s(x), \\
    (\nabla_X s)(x) &= \nabla_{X(x)} s.
  \end{align*}
  Also given a map \(g : Q \to N\) and \(Y_y \in T_y Q\) we have:
  \[
  \nabla_{d g Y_y} s = \nabla_{Y_y} (s \circ g).
  \]
  Moreover, if \(Y \in \mf{X}(Q)\) and \(X \in \mf{X}(N)\) are
  \(g\)--related in that \(d g(Y)(y) = X(g(y))\) for all \(y \in Q\)
  then:
  \[
  \nabla_Y (s \circ g) = (\nabla_X s) \circ g.
  \]

  We now turn to the covariant derivative on \(T L M\). 

  Let \(\alpha \in L M\) and let \(\beta, \gamma \in T_\alpha L M\) be
  never-zero and never-colinear.  Let \(t \in S^1\) and put \(x :=
  \alpha(t) \in M\).  Let \(\eta : T M \to M\) be a local addition on
  \(M\).  It is straightforward to show that the restriction of
  \(\eta\) to \(T_x M\) is a diffeomorphism onto an open neighbourhood
  of \(x\) of \(M\).  Write this as \(\eta_x : T_x M \to V_x\).

  As \(\alpha^* T M\) is orientable, we can choose a trivialisation
  \(\alpha^* T M \cong S^1 \times T_x M\).  From our assumptions on
  \(\beta\) and \(\gamma\) we can arrange matters so that under the
  induced isomorphism \(\Gamma_{S^1}(\alpha^* T M) \cong L T_x M\),
  they are taken to constant loops.  Define \(\psi : V \to L M\) as:
  \[
  V \xrightarrow{{\eta_x}^{-1}} T_x M \subseteq L T_x M \cong
  \Gamma_{S^1}(\alpha^* T M) \xrightarrow{\Psi_\alpha} U_\alpha
  \subseteq L M.
  \]
  Here, \(T_x M \subseteq L T_x M\) is identified with the subspace of
  constant loops.

  To determine the image of \(x\), we see that \({\eta_x}^{-1}(x) = 0
  \in T_x M\) which gets mapped to the zero loop and hence to the zero
  section, which is then mapped to \(\alpha\).  Hence \(\psi(x) =
  \alpha\).

  Apart from the map \(T_x M \to L T_x M\), all the maps are
  diffeomorphisms.  Thus to identify the image of \(d \psi_x\) it is
  sufficient to identify the image of the derivative of \(T_x M \to L
  T_x M\) at the origin.  As this map is linear, its derivative is
  itself and so we get the constant loops.  Thus, by construction, the
  image of \(d \psi_x\) contains \(\beta\) and \(\gamma\).  Let \(u, v
  \in T_x M\) be their respective pre-images.

  Using the charts \(\eta_x : T_x M \to V\) and \(\Psi_\alpha :
  \Gamma_{S^1}(\alpha^* T M) \to U_\alpha\) we can choose vector
  fields \(X_\beta\), \(X_\gamma\) on \(L M\) extending \(\beta\) and
  \(\gamma\) and \(Y_u\), \(Y_v\) on \(M\) extending \(u\) and \(v\)
  that are \(\psi\)--related and such that \([X_\beta, X_\gamma] =
  0\), \([Y_u, Y_v] = 0\).  Thus, writing \(\nabla^L\) for \(\nabla^{T
  L M}\):
  \begin{align*}
    \tau(\beta, \gamma) &= (\nabla^L_{X_\beta} X_\gamma)(\alpha) -
    (\nabla^L_{X_\gamma} X_\beta)(\alpha) \\
    &= (\nabla^L_{X_\beta} X_\gamma)(\psi(x)) - (\nabla^L_{X_\gamma}
    X_\beta)(\psi(x)) \\
    &= \big(\nabla^L_{Y_u} (X_\gamma \circ \psi)\big)(x) -
    \big(\nabla^L_{Y_v} (X_\beta \circ \psi)\big)(x).
  \end{align*}

  As \(\tau\) is a map into \(T_\alpha L M \subseteq T L M \cong L T
  M\), we can take its adjoint, \(\tau^\lor\).  Taking adjoints on the
  right, we use theorem~\ref{th:covadj} to simplify as follows:
  \[
  \big(\nabla^L_{Y_u} (X_\gamma \circ \psi)\big)^\lor =
  \nabla_{Y_u}(X_\gamma \circ \psi)^\lor.
  \]
  Now \(X_\gamma \circ \psi : M \to L M \to T L M\) adjoints to:
  \[
  S^1 \times M \xrightarrow{1 \times \psi} S^1 \times L M
  \xrightarrow{1 \times X_\gamma} S^1 \times L T M \xrightarrow{e} T
  M.
  \]
  Hence \((X_\gamma \circ \psi)^\lor(t, x) = e(t, X_\gamma(\alpha))
  = \gamma(t) = v\).  Thus we can consider \((X_\gamma \circ
  \psi)^\lor\) as a vector field on \(S^1 \times M\) extending \(v\)
  at \((t,x)\).  We extend the covariant derivative to \(S^1 \times
  M\) using the product of the connection on \(M\) with the standard
  connection on \(S^1\).  As the \(S^1\)--part is torsion free, the
  torsion of this extension is the torsion of \(M\).  Thus:
  \[
  \nabla_u (X_\gamma \circ \psi)^\lor - \nabla_v (X_\beta \circ
  \psi)^\lor = \tau(u,v).
  \]
  Hence \(\tau(\beta, \gamma)(t) = \tau(\beta(t), \gamma(t))\) and
  thus the torsion of the looped connection is the loop of the torsion
  of the original connection.
\end{proof}

\begin{corollary}
  \label{cor:torsion}
  A torsion-free connection on \(M\) loops to a torsion-free
  connection on \(L M\).
\end{corollary}

\subsection{The Enemy of my Enemy is not my Friend}
\label{sec:dual}

Section~\ref{sec:vect} showed that the loop space of a vector bundle
has a pleasant \(L \R\)-module structure.  We extend this further by
observing that:
\[
\hom_{L \R} (L E, L \R) = L \hom_\R(E,\R),
\]
and more generally:
\begin{align*}
\hom_{L \R} (L E, L F) &= L \hom_\R(E, F), \\
\iso_{L \R} (L E, L F) &= L \iso_\R(E,F).
\end{align*}
This last is a generalisation of the fact that if \(Q\) is the frame
bundle of \(E\) then \(L Q\) is the \(L \R\)-frame bundle of \(L E\).

However, when viewing \(L E\) as a mere vector bundle this pleasant
functorality is not preserved.  The most distressing case of this is
that the cotangent bundle of the loop space is not the loop space of
the cotangent bundle.  The simple explanation for this is that \(L
(E^*)\) is modelled on the space \(L \R^n\), as is \(L E\), but \((L
E)^*\) is modelled on its dual which is the space of \(\R^n\)-valued
distributions on the circle%
\footnote{%
  We are fortunate in this situation that dualising the dual -- which
  results in the \emph{bidual} -- ends us up where we started.  This
  is \emph{not} always true in infinite dimensions.
}.

Thus although \(T L M = L T M\) it is not true that \(T^* L M = L T^*
M\).  This rules out any possibility of an isomorphism between the
tangent and cotangent bundles, either via an inner product or via a
symplectic structure.  One can (easily) define \emph{weak} such
structures where the induced map \(T L M \to T^* L M\) is injective
but it can never be an isomorphism.

A choice of linear map \(f : L \R \to \R\) defines a map \(L (E^*) \to (L
E)^*\) via:
\[
L (E^*) = \hom_{L \R}(L E, L \R) \subseteq \hom_\R (L E, L \R)
\xrightarrow{f} \hom_\R(L E, \R) = (L E)^*.
\]
The most usual choice is the map \(\alpha \to \int_{S^1} \alpha\).
This has various properties that bode well for further construction,
not least being its equivariance under the natural action of the
circle.  With this choice, the induced map \(L (E^*) \to (L E)^*\) is
injective.  Using this, any isomorphism \(E \to E^*\) loops to give an
injection \(L E \to (L E)^*\), but never an isomorphism.

Thus over a loop space we have a schizophrenia as to whether bundles
should be vector spaces or \(L \R\)-modules.  On the one hand, most of
the constructions that one wishes to generalise from finite
dimensional topology definitely use vector spaces and cannot be
modified to use \(L \R\)-modules -- for the simple reason that the
differential of an arbitrary function \(f : L M \to \R\) is an
\R-linear map \(T L M \to \R\) and not \(L \R\)-linear.  On the other
hand, the theory of bundles as \(L \R\)-modules is very nice.  In
addition to the properties outlined above there is also the fact that
one never wants to consider \emph{any old} vector bundle with infinite
dimensional fibres as the corresponding full general linear group is
usually either not a Lie group or is contractible, neither of which is
much help.  Therefore we usually work with a subgroup, such as \(L
\gl(\R^n)\), which is a Lie group, does have some interesting topology,
and coincidentally preserves the \(L \R\)-module structure.

\newpage

\section{Submanifolds and Tubular Neighbourhoods}
\label{sec:subtube}

There are two important sources of submanifolds of loop spaces: those
that arise from coincidences and those that arise from the obvious
circle action.  The ``coincidental'' submanifolds are ones where some
constraint is imposed on the value of a loop at some specific times;
the most obvious being the based loop space.  The submanifolds arising
from the circle action are the fixed point sets of the various
subgroups of the circle.  The main result of this section is that all
of these submanifolds have tubular neighbourhoods.  We conclude this
section by exhibiting a submanifold without a tubular neighbourhood.

\subsection{The Fundamental Fibration}
\label{sec:funfibre}

The space of \emph{based} loops is one of the key concepts in
algebraic topology.  All of the above analysis works equally well for
the based loop space as for the free loop space.  The model spaces for
\(\Omega M\) are sections of \(\alpha^* T M\) which are zero at time
\(0\).  The tangent space of \(\Omega M\) is thus \(\Omega T M\) where
the base point in \(T M\) is the zero above the base point in \(M\).
The inclusion of based loops in free loops is smooth and fits into the
sequence:
\[
\Omega M \to L M \xrightarrow{e_0} M
\]
where \(e_0\) is the map which evaluates a loop at time \(0\).  This
is a fibration sequence and is split since \(M\) includes in \(L M\)
as the subspace of constant loops.

This is great as far as algebraic topology is concerned.  However for
differential topology we would like to know that this is a locally
trivial fibration.

\begin{theorem}
  \label{th:basedlt}
  Let \(M\) be a connected smooth manifold with base point
  \(x_0\).  Then \(\Omega M \to L M \to M\) is a locally trivial
  fibration.
\end{theorem}

In the non-connected case we simply use this result together with the
fact that the loop spaces -- based or free -- of a disjoint union are
a disjoint union of the loop spaces of the components.  Thus \(L M \to
M\) remains a locally trivial fibration although the fibre may differ
on components.

\begin{proof}
  We need to show local triviality.  This will imply that the
  diffeomorphism type of the fibres is locally constant and thus will
  allow us to identify every fibre with \(\Omega M\).

  We start with a technical result on diffeomorphisms of \(\R^n\).
  What we are looking for is a family \(\{\psi_v\}\) of compactly
  supported diffeomorphisms of \(\R^n\) indexed by points of \(\R^n\)
  such that \(\psi_v(0) = v\).  In other words, we are looking for a
  splitting of the map \(\Diff_c(\R^n) \to \R^n\) which evaluates a
  diffeomorphism at the origin.

  We construct this using the exponential map for compactly supported
  vector fields.  It is a standard corollary of ODE theory that the
  exponential map \(\exp : \m{X}_c(\R^n) \to \Diff_c(\R^n)\) is
  well-defined, where \(\m{X}_c(\R^n)\) is the space of compactly
  supported vector fields on \(\R^n\).

  We define a map \(\R^n \to \m{X}_c(\R^n)\) as follows: let \(\rho :
  \R \to [0,1]\) be a bump function such that:
  \[
  \rho(t) = \begin{cases}
    1 & 0 \le t \le 1 \\
    0 & 2 \le t.
	    \end{cases}
  \]
  Define \(\R^n \to \m{X}_c(\R^n)\), \(v \to X_v\) by:
  \[
  X_v(u) = \rho(\norm[u]^2)v.
  \]

  The exponential map \(\exp : \m{X}_c(\R^n) \to \Diff_c(\R^n)\) is
  such that the path \(t \to \exp(tX)y_0\) is the solution to the ODE
  \(y' = X(y)\) with initial condition \(y_0\).  Thus \(\exp(X_v)0\)
  is the value at time \(1\) of the solution to the ODE \(y' = v\)
  with initial condition \(y_0 = 0\).  Hence \(\exp(X_v)0 = v\).  Thus
  the map \(v \to \exp(X_v)\) is the required splitting.  Note that
  \(X_0\) is the zero vector field and so \(\exp(X_0)\) is the
  identity.

  Let \(x \in M\) and let \(\Omega_x M\) be the fibre of \(L M \to M\)
  at \(x\).  Let \(\phi : \R^n \to U\) be a chart at \(x\) with
  \(\phi(0) = x\).  Now \(\phi\) takes any compactly supported
  diffeomorphism of \(\R^n\) to one of \(U\).  As such a
  diffeomorphism is compactly supported, it is the identity near the
  boundary of \(U\).  It therefore extends to a diffeomorphism of
  \(M\) by defining it to be the identity outside \(U\).

  Using the above family of diffeomorphisms of \(\R^n\) we thus have a
  smooth family of diffeomorphisms of \(M\) indexed by the points of
  \(U\) with the property that \(\phi_u(x) = u\).
  
  Let \(L_U M := \{\alpha \in L M : \alpha(0) \in U\}\).  This is an
  open submanifold of \(L M\).  Define \(\Omega_x M \times U \to L_U
  M\) by \((\alpha, u) \to \phi_u(\alpha)\).  Since \(\phi_u(x) = u\),
  this is a fibrewise map of spaces over \(U\).  Its inverse is
  \(\beta \to (\phi_{\beta(0)}^{-1}(\beta), \beta(0))\).  This is the
  required local trivialisation.
\end{proof}

The key part of this proof is showing that the submanifold \(\Omega
M\) of \(L M\) has a tubular neighbourhood.  The length and intricacy
of this proof should be contrasted with the corresponding statement
about the submanifold \(M\) of \(L M\), where \(M\) is identified with
the space of constant loops in \(L M\).

\begin{proposition}
  The inclusion \(M \to L M\) admits a tubular neighbourhood.  The
  normal bundle is \(\Omega^v T M\), the space of fibrewise loops in
  \(T M\) which are based at the zero section; that is, the fibre of
  \(\Omega^v T M\) at a point \(p\) is \(\Omega T_p M\).  This can be
  identified with \(T M \otimes \Omega \R\).
\end{proposition}

\begin{proof}
  Let \(\eta : T M \to M\) be a local addition on \(M\) and let \(V
  \subseteq M \times M\) be the corresponding neighbourhood of the
  diagonal.  Since \(\Omega^v T M\) is a subset of the set of smooth
  loops in \(T M\), we can define \(\eta : \Omega^v T M \to L M\) by
  composition.  This is clearly a smooth map.

  Let \(p \in M\).  The domain of the chart map defined by \(\eta\) at
  the constant map at \(p\) can be naturally identified with \(L T_p
  M\).  The map \(\eta : \Omega^v T M \to L M\) restricted to the
  fibre above \(p\) is the restriction of the chart map to the
  subspace \(\Omega T_p M\).  Therefore \(\eta : \Omega^v T M \to L
  M\) is injective when restricted to any fibre.  The images of the
  fibres can be distinguished in \(L M\) since for \(\alpha\) in the
  fibre of \(\Omega^v T M\) above \(p \in M\), \(\eta \alpha(0) =
  \eta(O_p) = p\).

  The image of this map is the set of \(\alpha \in L M\) such that
  \((\alpha(0), \alpha(t)) \in V\) for all \(t \in S^1\).  This is
  open in \(L M\) as it is the preimage of \(L V\) under the
  continuous map \(L M \to M \times L M \to L M \times L M\) given by
  sending \(\alpha\) to \((\alpha(0), \alpha)\).

  The inverse of this map is thus \(\alpha \to (\pi \times
  \eta)^{-1}(\alpha(0), \alpha)\).  It is therefore a diffeomorphism
  onto its image.
\end{proof}

This is not the tubular neighbourhood of \(M\) that is usually
wanted as it is not \(S^1\)-equivariant.  We postpone the construction
of that neighbourhood to section~\ref{sec:cirtub}.

\subsection{Tubular Neighbourhoods}

The vector fields that we defined in section~\ref{sec:funfibre} on
\(\R^n\) did not use any structure of \(\R^n\) beyond its being an
inner product space.  Therefore we can define similar vector fields on
a vector bundle over a manifold.  We can use this to prove a
generalisation of this result involving tubular neighbourhoods.  One
important application of this generalisation is the following result:

\begin{proposition}
  \label{prop:diagloop}
  Let \(L M \times_M L M\) be the family of pairs of loops which
  coincide at time \(0\).  Then \(L M \times_M L M \to L M \times L
  M\) is an embedded submanifold with a tubular neighbourhood.
\end{proposition}

This result is used in \cite{rcjj} in the construction of the loop
product in the cohomology of the loop space.  A generalisation of it
is used in \cite{rcvg} to defined the other operations of
string topology in the cohomological setting.

We shall prove a little more than that a tubular neighbourhood
exists.  We shall prove that the \emph{obvious} neighbourhood is a
tubular neighbourhood.  To explain this remark, observe that there is
a pull-back diagram:
\[
\begin{CD}
  L M \times_M L M @>>> L M \times L M \\
  @Ve_0VV @VVe_0 \times e_0 V \\
  M @>\Delta>> M \times M
\end{CD}
\]
where \(\Delta\) is the inclusion of the diagonal.  The lower
line is an embedded submanifold with a tubular neighbourhood, say
\(V\), so define:
\[
L M \times_V L M : \{(\alpha, \beta) \in L M \times L M : (\alpha(0),
\beta(0)) \in V\}.
\]
This fits in to the above diagram very neatly:
\[
\begin{CD}
  L M \times_M L M @>>> L M \times_V L M \subseteq L M \times L M \\
  @Ve_0VV @VVe_0 \times e_0 V \\
  M @>\Delta>> V \subseteq M \times M
\end{CD}
\]
It would be nice if not only did \(L M \times_M L M\) have a tubular
neighbourhood in \(L M \times L M\) but that \(L M \times_V L M\) were
an example of such.

Before proving that this is so, let us examine what we get for free
and thus what extra is needed to be shown.  As the normal bundle to
the diagonal embedding is isomorphic to \(T M\), the existence of the
lower tubular neighbourhood means that there is a diffeomorphism \(T M
\to V \subseteq M \times M\) where \(V\) is an open neighbourhood of
the diagonal such that the composition of this with the zero section
is the embedding of the diagonal.  For convenience, we shall assume
that this diffeomorphism comes from a local addition on \(M\), see
definition~\ref{def:locadd}.  Thus we have an identification between
tangent vectors and pairs of suitably close points.  Our assumption
means that the first of those points is the anchor for the tangent
vector.

What we do get for free is that the normal bundle on the upper level
is the pull-back of the normal bundle on the lower level.  Thus on the
upper level we seek an identification between the spaces:
\begin{align*}
{e_0}^*T M &= \{(\alpha, \beta, v) : \alpha(0) = \beta(0), v \in
T_{\alpha(0)} M\}, \\
L M \times_V L M &= \{(\alpha, \beta) : (\alpha(0), \beta(0)) \in V\}.
\end{align*}
To make everything fit nicely into the pull-back diagram, we want the
diffeomorphism between these two to project down to the diffeomorphism
that we already have.   Since \(\alpha(0)\) therefore should not
change, we may as well -- for simplicity -- assume that \(\alpha\)
does not move.  Thus we want \((\alpha, \beta, v) \to (\alpha,
\tilde{\beta})\) such that on evaluation at \(0\) we get the lower
diffeomorphism.

The difficulty is that \(v\) only tells us what to do with
\(\beta(0)\).  As \(\beta\) is smooth, we need to know what to do with
the rest of it.  This involves some choices and some careful analysis.
Fortunately, we have already laid the necessary foundations.

\begin{proposition}
  Let \(M\) be a smooth finite dimensional manifold, \(P \subseteq M\)
  an embedded submanifold with normal bundle \(E\) and tubular
  neighbourhood \(V \subseteq M\) with diffeomorphism \(\nu : E \to
  V\).  Let \(L_P M := \{\alpha \in L M : \alpha(0) \in P\}\) and
  \(L_V M := \{\alpha \in L M : \alpha(0) \in V\}\).

  The inclusion \(L_P M \to L M\) is a smooth embedding with normal
  bundle \({e_0}^* E\) and tubular neighbourhood \(L_V M\).  Moreover,
  there is a diffeomorphism \({e_0}^* E \to L_V M\) covering \(\nu : E
  \to V\).
\end{proposition}

We view \(P\) as an actual subset of \(M\) rather than taking \(j : P
\to M\) as an embedding to reduce the number of maps that we need to
make explicit.

\begin{proof}
  We omit the full proof that \(L_P M\) is a submanifold of \(L M\).
  The proof that it is a manifold is a repetition of the proof that
  \(L M\) is a manifold.  The embedding follows from the fact that the
  isomorphism: \(L \R^k \cong \Omega \R^k \oplus \R^k\) together with
  the fact that \(P\) is embedded in \(M\).  The case of \(L_V M\) is
  simpler as it is an open subset of \(L M\).

  The bundle \({e_0}^* E \to L_P M\) is the pull-back bundle via the
  evaluation map \(\alpha \to \alpha(0)\).  As a space,
  \[
  {e_0}^* E = \{(\alpha, v) \in L_P M \times E : \alpha(0) = \pi(v)\}.
  \]

  We equip \(E\) with inner products on the fibres, varying smoothly
  over \(P\).  The pull-back, \(e_0^* E\) inherits these inner
  products.  Let \(\norm\) be the corresponding fibrewise norm.

  Using the local triviality of \(E\) and paracompactness of \(P\) we
  wish to choose an open cover of \(P\) over which \(E\) trivialises
  together with a variation on the theme of a subordinate partition of
  unity.  The variation that we want is that the \emph{squares} of our
  functions should be a partition of unity.  This presents no
  technical difficulties: recall that the final step in constructing a
  partition of unity is to renormalise a family of bump functions with
  respect to their sum; if one instead renormalised with respect to
  the square-root of the sum of their squares, the resulting family
  would have the required property.  This square-root results in a
  smooth function as it is the square-root of a strictly positive
  function.

  Thus we choose, for an indexing set \(\Lambda\):
  \begin{enumerate}
  \item an open cover \(\{U_\lambda : \lambda \in \Lambda\}\),
  \item trivialisations \(\phi_\lambda : E_\lambda \to U_\lambda
    \times \R^k\),
  \item smooth functions \(\rho_\lambda : P \to \R\) with compact
    support such that \(\{{\rho_\lambda}^2\}\) is a partition of unity
    subbordinate to \(\{U_\lambda\}\) with the support of
    \(\rho_\lambda\) contained in \(U_\lambda\).
  \end{enumerate}

  Let \(\tilde{\phi}_\lambda : E_\lambda \to \R^k\) be the
  composition of \(\phi_\lambda\) with the projection onto \(\R^k\).
  Define \(s : E \to \Gamma(E)\) by:
  \[
  s(v)(x) = \sum_{\lambda \in \Lambda} \rho_\lambda(\pi(v))
  \rho_\lambda(x) \phi_\lambda^{-1}(x, \tilde{\phi}_\lambda(v)).
  \]
  Note that the \(\lambda\)-summand is zero unless both
  \(\rho_\lambda(x)\) and \(\rho_\lambda(\pi(v))\) are non-zero.
  Therefore the support of the section \(s(v)\) is contained in the
  union of the supports of the \(\rho_\lambda\) for which
  \(\rho_\lambda(\pi(v)) \ne 0\).  As the supports of the
  \(\rho_\lambda\) form a locally finite family of compact sets, the
  support of \(s(v)\) is compact and hence \(s\) takes values in
  \(\Gamma_c(E)\), sections with compact support.

  This function is smooth and has the following properties:
  \begin{enumerate}
  \item the restriction to a fibre is linear, and
  \item \(s(v)(\pi(v)) = v\).
  \end{enumerate}
  The first of these follows from the fact that the \(\phi_\lambda\)
  are linear on fibres.  The second uses the fact that the
  \(\rho_\lambda\) square to a partition of unity.  Note that as a
  consequence we have that if \(v\) is a zero vector then \(s(v)\) is
  the zero section.

  From this family of sections of \(E\), we define a family of
  compactly supported vector fields on \(E\) which, on fibres, look
  like the vector fields that we used in the proof that \(L M \to M\)
  was locally trivial.  We note that the tangent bundle of \(E\)
  contains a canonical copy of \(\pi^* E\) as the vertical tangent
  bundle.  Thus a section \(\sigma\) of \(E\) defines a vector field
  on \(E\) by \(v \to \sigma(\pi(v))\).  Let \(\tau : \R \to [0,1]\)
  be a bump function with \(\tau(t) = 1\) for \(0 \le t \le 1\) and
  \(\tau(t) = 0\) for \(t \ge 2\).  Define \(X : \Gamma(E) \to
  \m{X}(E)\) by:
  \[
  X_\sigma(v) = \tau(\norm[v]^2)\sigma(\pi(v)).
  \]
  By construction, \(X_s\) has fibrewise compact support and its
  horizontal support agrees with that of \(\sigma\).  Therefore this
  restricts to a map \(X : \Gamma_c(E) \to \m{X}_c(E)\).  Combining
  this with the above map \(E \to \Gamma_c(E)\) and the exponential
  map \(\Gamma_c(E) \to \Diff_c(E)\) we obtain a map \(\psi : E \to
  \Diff_c(E)\).

  Now the vector field corresponding to a point \(v \in E\) takes
  values in the vertical tangent space of \(E\).  Therefore the
  diffeomorphism is a fibre-preserving diffeomorphism covering the
  identity on the base.  On fibres, it looks like the diffeomorphisms
  we had in the previous proof.  Thus as \(s_v(\pi(v)) = v\),
  \(\psi(v)(0_{\pi(v)}) = v\), where \(0_x\) is the zero vector in
  \(E_x\).

  The diffeomorphism \(\nu : E \to V\) defines \(\Diff_c(E) \to
  \Diff_c(V)\) and thence to \(\Diff_c(M)\) since a compactly
  supported diffeomorphism in \(V\) extends by the identity to the
  whole of \(M\).  Hence we have \(\theta : V \to \Diff_c(M)\) such
  that \(\theta(v)(\pi(\nu^{-1}v)) = v\).

  We now define our tubular neighbourhood diffeomorphism as:
  \[
  (\alpha,v) \to \theta(\nu(v))(\alpha)
  \]
  with inverse:
  \[
  \beta \to (\theta(\beta(0))^{-1}(\beta), \nu^{-1}(\beta(0))).
  \]
  Evaluating at zero, on the left we get \(v \in E_{\alpha(0)}\)
  whilst on the right we get \(\theta(\nu(v))(\alpha(0))\).  Now as
  \(\pi(v) = \alpha(0)\), \(\theta(\nu(v))(\alpha(0)) =
  \theta(\nu(v))(\pi(v)) = \nu(v)\).  Hence the diffeomorphism on the
  loop spaces projects down to \(\nu : E \to V\) under evaluation. 
\end{proof}

Proposition~\ref{prop:diagloop} follows immediately using \(\Delta : M
\to M \times M\).  Theorem~\ref{th:basedlt} also follows from this
proposition using the embedding of a single point.

\subsection{Equivariant Tubular Neighbourhoods}
\label{sec:cirtub}

The previous section deals with submanifolds arising from
``coincidences'': loops that happen to coincide with each other or
with some submanifold of the target manifold.  Another source of
submanifolds comes via the natural circle action on the loop space.

\begin{defn}
  Define the circle action \(\rho : S^1 \times L M \to L M\) by
  \(\rho(t, \alpha)(s) = \alpha(t + s)\).
\end{defn}

The adjoint of \(\rho\) is the composition:
\[
S^1 \times S^1 \times L M \xrightarrow{(s,t,\gamma) \to (s + t,
  \gamma)} S^1 \times L M \xrightarrow{e} M
\]
which is obviously smooth.

This action induces an action by any subgroup of \(S^1\).  We shall be
concerned with the \emph{compact} subgroups which are the finite
cyclic groups and \(S^1\) itself.  We wish to consider the fixed point
subsets of these actions.  It is straightforward to show that for \(G
= S^1\) the fixed points are the constant loops and so the fixed point
set is diffeomorphic to \(M\) while for \(G \ne S^1\) the fixed point
set is the set of loops of period \(1/\abs{G}\) and this is
diffeomorphic to \(L M\).  What is more intricate is showing that
these all have \(S^1\)-equivariant tubular neighbourhoods.

\begin{theorem}
  \label{th:eqtub}
  Let \(G \subseteq S^1\) be a compact subgroup (including the case
  \(G = S^1\)).  The fixed point set of the induced action of \(G\) on
  \(L M\) is an \(S^1\)-invariant embedded submanifold with an
  \(S^1\)-equivariant tubular neighbourhood.  The normal bundle is \(T
  (L M^G) \otimes \Ci(G, \R)_0\) where \(\Ci(G, \R)_0\) is the
  \(G\)-invariant complement of the constant maps in \(\Ci(G, \R)\).
\end{theorem}

Strangely, for \(G \ne S^1\) the results about the \(G\)-fixed points
in \(L M\) depend on the structure of \(\map(G, M)\).  For more on the
links between loop spaces and the spaces \(\map(G, M)\) see~\cite{jj}
and~\cite{as3}.  Note that as \(G\) is finite, \(\map(G, M)\) is a
product of copies of \(M\).

To prove this theorem we need more structure on \(M\).  In studying
the whole loop space, \(L M\), we used a \emph{local addition} to
enable us to use the tangent spaces.  Now we need a \emph{local
averaging function}.  The point is that we need to be able to find an
\(S^1\)-equivariant map from sufficiently small loops in \(M\) to
\(M\).  When considering the non-equivariant tubular neighbourhood of
the constant loops in \(L M\) we could just take evaluation at a
point.  This is not \(S^1\)-equivariant so is not adequate for our
purposes.  In \(\R^n\) we would average the values taken by the loop
(equivalently, take the constant Fourier component).  A \emph{local
averaging function} is precisely what we need in order to extend this
to an arbitrary manifold.  It is somewhat more complicated than a
local addition and so we shall give an explicit construction rather
than a definition.

The starting point for this construction is an embedding of \(M\) in
some Euclidean space \(\R^k\).  We shall identify \(M\) with its image
to avoid excess maps.  This also identifies \(T M\) with its image in
\(T \R^k\).  We consider \(T_x \R^k\) to be an affine space anchored
at \(x\) and isomorphic to \(\R^k\).  Thus the addition in \(T_x
\R^k\) is \((u,v) \to (u - x) + (v - x) + x = u + v - x\).

This identifies \(T_p M\) with an affine subspace of \(\R^k\) anchored
at \(p \in M\).  Let \(\pi : N \to M\) be the vector bundle defined by
setting \(N_p\) to be the affine orthogonal complement to \(T_p M\)
(also anchored at \(p\)).  Thus as an affine space, \(N_p = p + (T_p M
- p)^\perp\).

For \(p \in M\) we have an orthogonal projection map \(\R^k \to T_p M\)
which we restrict to \(M\) to define \(\lambda_p : M \to T_p M\).
This map varies smoothly in \(p\) so we define \(\lambda : M \times M
\to T M\) by:
\[
\lambda (p,q) = \lambda_p(q) \in T_p M.
\]

This map has the following properties:
\begin{enumerate}
\item Consider \(M \times M\) as a bundle over \(M\) via projection
  onto the first factor.  Then \(\lambda\) is a bundle map.

\item The composition of \(\lambda\) with the diagonal map
  \(M \to M \times M\) is the zero section of \(T M\).
  \label{it:zero}

\item The derivative at \((p,p) \in M \times M\) is an isomorphism.
  This is because \(d \lambda_p\) at \(p\) is the identity.
\end{enumerate}

Standard techniques of differential topology involving judicious use
of the inverse function theorem thus allow us to find a neighbourhood
of the diagonal in \(M \times M\) on which \(\lambda\) restricts to a
diffeomorphism.  We therefore have \(\lambda : M \times M \supseteq V
\xrightarrow{\cong} U \subseteq T M\).  Let \(\eta : T M \supseteq U
\to M\) be the composition of \(\lambda^{-1}\) with the projection
onto the second factor.  As \(\lambda\) is a bundle map, the
projection onto the first factor is just \(\pi : T M \to M\) so \(\pi
\times \eta = \lambda^{-1}\).  By item~\eqref{it:zero} above, \(\eta\)
composed with the zero section is the identity on \(M\).  Thus
\(\eta\) is almost a local addition, the only variation is that its
domain is not the whole of \(T M\).

There is a natural map \(\iota : N \to \R^k\) given by the natural
inclusions \(N_p \to \R^k\).  The derivative of \(\iota\) at a point
in the image of the zero section is an isomorphism as it corresponds
to the isomorphism \(T_p M \oplus N_p = T_p \R^k\).  The composition
of \(\iota\) with the zero section \(M \to N\) is just the embedding
\(M \to \R^k\).  Thus by similar techniques of differential topology,
there is a neighbourhood \(W\) of the zero section in \(N\) and a
neighbourhood \(X\) of \(M \subseteq \R^k\) such that \(\iota : W \to
X\) is a diffeomorphism.  Thus we have a map \(\pi \iota^{-1} : X \to
M\) with the property that for \(x \in X\), \(x - \pi \iota^{-1}(x)\)
is orthogonal to \(T_{\pi\iota^{-1}(x)} M\).

The two maps \(\eta\) and \(\iota\) have been defined without
reference to each other.  We shall use them together so we need to
modify their domains and codomains so that they interact nicely.  The
modification that we need to make is to shrink \(U\) so that the
closure of the convex hull of \(\eta(U \cap T_p M)\) -- taken in
\(\R^k\) -- is contained in \(X\), the codomain of \(\iota\).  This
ensures that if \(C \subseteq U \cap T_p M\) is any set then the
closure of the convex hull of \(\eta(C)\) is contained in the codomain
of \(\iota\).

We shall now explain how we are going to use these two maps to
construct the required tubular neighbourhoods.  Let \(G\) be a compact
subgroup of \(S^1\).  Let \(\Ci_0(G, U)\) be the space of smooth maps
\(\alpha : G \to U\) with the following properties:
\begin{enumerate}
\item There is some \(p \in M\) depending on \(\alpha\) such that
  \(\alpha(G) \in U \cap T_p M\).

\item Considering \(\alpha\) as a map into \(\R^k\) via the inclusion
  \(T_p M \to \R^k\), the following holds:
  \[
  \int_G \alpha = p.
  \]
\end{enumerate}
In this last property, if \(G\) is finite then the integral is simply
the average value of \(\alpha\) on \(G\), if \(G = S^1\) then we take
the standard \(S^1\)-invariant measure of total volume \(1\).  Note
that in our view \(p\) is the zero in \(T_p M\).

Composition with \(\eta\) defines \(\Ci_0(G, U) \to \Ci(G, M)\).

\begin{lemma}
  The map \(\Ci_0(G,U) \to \Ci(G,M)\) is a diffeomorphism onto its
  image which is open in \(\Ci(G,M)\).
\end{lemma}

\begin{proof}
  Let \(\cvx(G,M)\) denote the family of smooth maps \(\beta : G \to
  M\) such that when considered as a map into \(X\) the closed convex
  hull of \(\beta(G)\) lies inside \(X\).

  As \(G\) is compact and \(\beta\) continuous, the closed convex hull
  of \(\beta(G)\) is compact.  Therefore \(\cvx(G,M)\) is open within
  \(\Ci(G,M)\).  For \(\beta \in \cvx(G,M)\), the value of \(\int_G
  \beta\) lies within the closed convex hull of \(\beta(G)\) and
  therefore in \(X\).  Hence we have a well-defined smooth map \(\tau
  : \cvx(G,M) \to M\) given by \(\tau(\beta) = \iota^{-1} \int_G
  \beta\).

  Consider the set:
  \[
  Y_G := \{\beta \in \cvx(G, M) : (\tau \beta, \beta(t)) \in V \text{
  for all } t \in S^1\}.
  \]
  We assume that \(V \subseteq M \times M\) was modified at the same
  time as \(U\) so that \(\lambda : U \to V\) remains a
  diffeomorphism.  The set \(Y_G\) is open in \(\cvx(G, M)\), whence
  in \(L M\), as it is the preimage of \(L V \subseteq L M \times L
  M\) under the map \(\cvx(G, M) \to M \times L M \to L M \times L
  M\), \(\beta \to (\tau \beta, \beta)\).

  Now for \(\beta \in Y_G\), \(\lambda(\tau \beta, \beta)\) is a
  map \(G \to U \subseteq T M\).  Since \(\pi \lambda(p,q) = p\),
  it takes values in the fibre \(U \cap T_{\tau \beta} M\).  As
  \(\lambda_p : M \to  T_p M\) is the orthogonal projection, the
  difference \(\beta - \lambda(\tau \beta, \beta)\) is orthogonal
  to \(T_{\tau \beta} M\).  Hence \(\int_G \lambda(\tau \beta,
  \beta) \in N_{\tau \beta}\).  Since \(T_{\tau \beta} M\) is convex,
  this integral also lies in \(T_{\tau \beta} M\).  It is therefore
  \(\tau \beta\), the zero point in \(T_{\tau \beta} M\).
  Hence \(\lambda(\tau \beta, \beta) \in \Ci_0(G, U)\).

  Now \(\eta \lambda(p,q) = q\) so \(\eta \lambda(\tau \beta, \beta)
  =\beta\).  Hence the map \(Y \to \Ci_0(G,U)\) is inverse to the map
  \(\alpha \to \eta \alpha\).  As these maps are both smooth, they are
  diffeomorphisms.
\end{proof}

\begin{lemma}
  There is a smooth \(G\)-equivariant map \(\Ci_0(G, T M) \to \Ci_0(G,
  U)\) which is a diffeomorphism onto an open subset.
\end{lemma}

\begin{proof}
  We build this map in two stages: fibrewise and then extend over
  \(M\).

  For the fibrewise situation, let \(I : \Ci(G, \R^n) \to \R^n\) be
  the integration map.  As this is \(G\)-equivariant with respect to
  the trivial \(G\) action on \(\R^n\), \(\ker I\) is a
  \(G\)-invariant subspace.  We wish to define a \(G\)-equivariant
  diffeomorphism \(\ker I \to \ker I \cap \Ci(G, D^n)\) where
  \(D^n\) is the open unit disc in \(\R^n\).

  We do this by noting that \(\Ci(G, D^n)\) is also the intersection
  of \(\Ci(G, \R^n)\) with the unit ball in \(C(G,\R^n)\), the space
  of continuous maps with the standard sup-norm.  A symmetric
  diffeomorphism \(\phi : \R \to (-1,1)\) thus defines a
  diffeomorphism \(\hat{\phi} : \Ci(G, \R^n) \to \Ci(G, D^n)\) via
  (for \(\alpha \ne 0\)):
  \[
  \alpha \to \frac{\phi(\norm[\alpha]_\infty)}{\norm[\alpha]_\infty}
  \alpha.
  \]

  Therefore \(I \hat{\phi}(\alpha) = 0\) if and only if \(I \alpha =
  0\) so \(\hat{\phi}\) preserves \(\ker I\).  It therefore defines a
  diffeomorphism \(\ker I \to \ker I \cap \Ci(G, D^n)\).

  We extend this over the manifold by choosing a smooth function
  \(\epsilon : M \to (0, \infty)\) such that the \(\epsilon\)-ball in
  \(T_p M\) is contained in \(U\).  The map \(\alpha \to \epsilon(p)
  \hat{\phi}(\alpha)\) defines the required diffeomorphism.
\end{proof}

The part of theorem~\ref{th:eqtub} for \(G = S^1\) follows immediately
by putting \(G = S^1\) in the above.  The rest of
theorem~\ref{th:eqtub} is equally simple but requires a word or two of
explanation.

The required neighbourhood of \(L M^G\) in \(L M\) consists of those
loops which, when evaluated on the cosets of \(G\) in \(S^1\), take
values in the neighbourhood \(Y_G\) of the constant maps in \(\Ci(G,
M)\).  We use the contraction of \(Y_G\) onto \(M\) to define the
corresponding contraction of the neighbourhood onto \(L M^G\).  Thus
by restricting a loop to each coset of \(G\) in turn we move from the
infinite case to the finite case.  The \(G\)-invariance of everything
in finite dimensions means that when evaluating on a coset you get the
same answer no matter which point you choose as the initial point,
thus the result is well-defined.

\subsection{A Not-So-Nice Submanifold}
\label{sec:notube}

We conclude this section with an example of a submanifold that does
not have a tubular neighbourhood.  As with so many counterexamples or
counterintuitive results in infinite dimensions, the failure is due to
a linear problem.

Let \(L_\flat M\) denote the space of loops in \(M\) that are
infinitely flat at the basepoint of \(S^1\).  We allow the value of
this basepoint to vary, the flatness condition is concerned with the
derivatives.  This is a smooth manifold modelled on \(L_\flat \R^n\)
and is a submanifold of \(L M\).  However, it does not posses a
tubular neighbourhood.

This is for the simple reason that the exact sequence:
\[
0 \to L_\flat \R \to L \R \to \R^\N \to 0
\]
does not split.  The map \(L \R \to \R^\N\) sends a map to its
derivatives at \(0\).  That this sequence is exact and does not split
is a corollary of~\cite[Lemma 21.5]{akpm}.  Therefore the inclusion
\(L_\flat M \to L M\) does not have a normal bundle.

\newpage

\section{A Miscellany}
\label{sec:misc}

We conclude this document with two topics designed to lead the
interested reader out of the basic differential topology of the loop
space and into more interesting areas.  The first topic is the
differential geometry of the loop space.  This is generally more
complicated than finite dimensional geometry, but still a certain
amount can be said without too much difficulty.  The second topic is
the semi-infinite structure of a loop space, which has no true analogy
in finite dimensions and is an important area of current interest.

\subsection{Weak Riemannian Manifolds}

Fairly early on in any text on Riemannian geometry is the statement
that any finite dimensional manifold admits a Riemannian structure.
We used this fact in demonstrating that any such manifold admits a
local addition.  In finite dimensions the definition of a Riemannian
structure is straightforward: it consists of a smooth choice of inner
product on each fibre of the tangent bundle.  In infinite dimensions
things are more complicated.  The following discussion clearly
generalises to inner products on arbitrary vector bundles.

The issues that one needs to deal with are:
\begin{enumerate}
\item Fibrewise questions:
  \begin{enumerate}
  \item Do the fibres admit (smooth) inner products?

    For example, the space of all \R-valued sequences with its inverse
    limit topology does not.

    For the model space of the loop space, \(L \R^n\), the answer is
    ``yes''.

  \item Up to equivalence, how many inner products are there?

    Equivalence means that there is a topological isomorphism taking
    one inner product to the other.  The answer is likely to be that
    standard mathematical answer: none, one, or infinity.

    For ``none'', the previous example works.  For ``one'' we takes
    its dual: the space of all \R-valued sequences that are eventually
    zero.  For ``infinity'' we can take the Hilbert space of
    square-integrable \R-values sequences.

    For the model space of the loop space, the answer is ``infinity''.

  \item Does a particular inner product induce an isomorphism to the
    dual space?

    The global version of this -- the induced isomorphism of the
    tangent and cotangent bundles -- is one of the mainstays of finite
    dimensional geometry since it allows free movement between vector
    fields and one-forms.  In infinite dimensions a positive answer to
    this question means that one is dealing with a Hilbert space (with
    the correct choice of inner product).  Therefore if one wishes to
    work with more than just Hilbert manifolds one must be prepared
    for a negative answer.

    What one \emph{always} has is a linear injection into the dual
    space.

    For the model space of the loop space, the answer is ``no''.
  \end{enumerate}

  Once one has answered these questions, and has a positive answer to
  the first, then the partition-of-unity argument applies and one can
  define a smooth global choice of inner product on the fibres of the
  tangent bundle.  However, the questions don't stop there:

\item Global questions:

  \begin{enumerate}
  \item Is this Riemannian structure \emph{weak} or \emph{strong}?

    The difference is whether or not the inner products identify the
    tangent and cotangent bundles, with ``strong'' meaning that they
    do.  A strong Riemannian structure is only possible when one has a
    Hilbert bundle and the inner product on each fibre is equivalent
    to the standard one.

    Thus with a weak structure all vector fields are one-forms but the
    converse only holds for a strong structure.

    For a loop space, the answer is always ``weak''.

  \item Is the equivalence class of the inner product (locally)
    constant?

    The problem here is that the partition of unity construction paid
    no attention to the question of equivalence.  As an example,
    consider a space \(E\) with two inequivalent inner products \(g_0\)
    and \(g_1\).  Let \(\rho : [0,1] \to [0,1]\) be the identity map,
    then \(\{\rho, 1 - \rho\}\) is a partition of unity on \([0,1]\).
    Define a fibrewise inner product on \([0,1] \times E\) by \(g(t)
    := \rho(t) g_1 + (1 - \rho(t)) g_0\).  By construction, the
    equivalence class of this inner product is not locally constant.

    This is closely related to the next question:

  \item Is there a bundle (i.e.~locally trivial) of Hilbert spaces
    which can be considered as the fibrewise completions of the fibres
    of the tangent bundle?

    The connection with the previous question comes about because an
    inner product defines a Hilbert completion.  An equivalence
    between two inner products extends to an isometric isomorphism of
    the corresponding completions.  Therefore if the equivalence class
    is locally constant the Hilbert completions fit together to define
    a locally trivial Hilbert bundle.

    Note that one can define this bundle without reference to an
    actual inner product but only to an equivalence class.
    Essentially, one breaks down the choice of inner product to an
    initial choice of equivalence class -- which defines the bundle of
    completions -- and then to a choice of inner product within that
    class.
    
  \item Is the tangent bundle with its family of inner products
    \emph{isometrically} locally trivial?

    By this we mean that there is one fixed inner product on the model
    space and the tangent bundle can be locally trivialised in such a
    way that the fibrewise inner products are all carried to this
    reference one.  In finite dimensions this follows from the
    Gram-Schmidt algorithm.

    A positive answer to this question implies a positive answer to
    the previous one since there is a fixed Hilbert completion
    corresponding to the fixed inner product.

  \item Assuming the existence of a smaller group than the full
    general linear group, can the construction be done in such a way
    that the transition functions lie in this group?

    If the answer to this question is yes (or is suspected to be),
    this provides a simpler route to the construction: first fix the
    reference inner product and, consequently, Hilbert completion.
    Then prove that the given group preserves this inner product, and
    hence the Hilbert completion.  Finally, use this group action to
    transfer the whole structure to the manifold.
  \end{enumerate}
\end{enumerate}

We illustrate this with the loop space, \(L M\), of a finite
dimensional Riemannian manifold \(M\).  There is a canonical weak
Riemannian structure on \(L M\) coming from the Riemannian structure
on \(M\).  There are two ways to define these inner products.

The direct way is to use the strategy of section~\ref{sec:dual}.  The
inner product on the tangent space of \(M\) is a symmetric fibrewise
bilinear map \(T M \times_M T M \to \R\) with the property that the
induced map \(g : T M \to T^* M\) satisfies \(g(v)(v) > 0\) for \(v
\ne 0\) (i.e.~not in the image of the zero section).  This loops to a
symmetric bi-\(L \R\)-linear map \(L T M \times_{L M} L T M \to L \R\)
such that the induced map \(g : L T M \to L T^* M\) satisfies
\(g(\alpha)(\alpha) > 0\) for \(\alpha \ne 0\).  The inequality now
holds in \(L \R\) and is defined by \(\beta > \gamma\) if \(\beta(t)
\ge \gamma(t)\) for all \(t\) and \(\beta \ne \gamma\) (equivalently,
there is some \(t\) such that the inequality is strict).  We then
apply the integration map \(\int_{S^1} : L \R \to \R\).  This is
an order-preserving linear map and so the symmetric bilinear map \(L T
M \times_{L M} L T M \to \R\) has the property that the induced map
\(\int g : L T M \to L T^* M \to T^* L M\) satisfies \(\int
g(\alpha)(\alpha) > 0\) for \(\alpha \ne 0\).  Untangling all of that
yields the formula:
\[
\ipv{\beta}{\gamma}_\alpha = \int_{S^1}
\ip{\beta(t)}{\gamma(t)}_{\alpha(t)} d t,
\]
for \(\beta, \gamma \in \Gamma_{S^1}(\alpha^* T M) = T_\alpha L M =
L_\alpha T M\).

The indirect way is to observe that the structure group of \(M\) is,
because of the choice of Riemannian structure, \(O_n\).  Therefore the
structure group of \(L M\) is \(L O_n\).  Now the action of \(L O_n\)
on \(L \R^n\) preserves the standard inner product coming from the
inclusion \(L \R^n \to L^2 \R^n\) (in fact, it is precisely the
subgroup of \(L \gl_n(\R)\) which does so).  Therefore we can define
a locally trivial inner product on the fibres \(T L M\) and a
corresponding bundle of Hilbert completions.

The equivalence of the two approaches comes from the fact that the
principal \(L O_n\)-bundle of \(T L M\) is the loop of the principal
\(O_n\)-bundle of \(T M\).  An element of the \(L O_n\)-bundle above
\(\alpha \in L M\) is an isometric trivialisation of the bundle
\(\alpha^* T M \to S^1\).  This defines an isometric isomorphism
\(\Gamma_{S^1}(\alpha^* T M) \to L \R^n\) (assuming orientability to
avoid twisting) and hence identifies the inner product given by the
above formula with the standard one.

We note that both approaches have their advantages.  In the first
there is an explicit formula for the inner product that one can work
with.  In the second, the local triviality and the existence of the
bundle of Hilbert completions are straightforward.

For this weak Riemannian structure on \(L M\), it is straightforward
to prove that certain geometric objects on \(M\) loop to the
corresponding objects on \(L M\).

\begin{proposition}
  The Levi-Civita connection on \(M\) loops to the Levi-Civita
  connection on \(L M\).
\end{proposition}

\begin{proof}
  From corollary~\ref{cor:torsion}, the loop of the Levi-Civita
  connection is torsion-free.  To see that it respects the inner
  product, we use the fact that the orthogonal structure group of \(T
  L M \cong L T M\) is the loop of the orthogonal structure group of
  \(T M\).  Hence as the Levi-Civita connection is an orthogonal
  connection, its loop is also orthogonal.
\end{proof}

The Koszul formula that is often used to prove existence and
uniqueness of the Levi-Civita connection can only, in infinite
dimensions, be used to prove uniqueness.  This is because the
existence part of the proof uses the isomorphism of vector fields and
one-forms at a crucial stage but the uniqueness only uses the
injectivity of the map from vector fields to one-forms.  Hence a
corollary of the above result is that the Levi-Civita connection on
\(L M\) does exist.

\begin{proposition}
  For \(X\) either \(M\) or \(L M\) and \(v \in T X\), let \(\gamma_v
  : I_v \to X\) denote the geodesic corresponding to \(v\) with
  maximal domain \(I_v\).  Then for \(\nu \in T L M\),
  \[
  e_t \gamma_\nu = \gamma_{e_t \nu},
  \]
  and \(I_\nu = \cap I_{e_t \nu}\).
\end{proposition}

\begin{proof}
  Let \(\gamma : I \to L M\) be a path.  For \(t \in S^1\), let
  \(\gamma_t : I \to M\) be the adjoint of \(\gamma\) restricted to
  \(\{t\} \times I\).  Differentiating, we get \(\gamma' : I \to T L
  M\) and its adjoint when restricted to \(\{t\} \times S^1\) is
  \({\gamma_t}'\).

  We have the covariant derivative of \(\gamma'\) along the canonical
  vector field \(\partial x\) of \(I\): \(\gamma'' :=
  \nabla^L_{\partial x} \gamma'\).  This is again a map \(I \to T L
  M\) above \(\gamma\).  By theorem~\ref{th:covadj}, this has adjoint:
  \[
  \big( \nabla^L_{\partial x} \gamma' \big)^\lor = \nabla_{\partial x}
  ({\gamma'}^\lor).
  \]
  Hence the adjoint of \(\gamma''\) restricted to \(\{t\} \times S^1\)
  is \({\gamma_t}''\).  Thus \(\gamma''\) vanishes if and only if each
  \({\gamma_t}''\) vanishes and so \(\gamma\) is a geodesic if and
  only \(e_t \gamma\) is a geodesic for each \(t \in S^1\).  The rest
  of the proposition then follows directly.
\end{proof}

\begin{corollary}
  If \(M\) is geodesically complete then \(L M\) is geodesically
  complete.
\end{corollary}

In fact, this is an ``if and only if'' as \(M\) is a Riemannian
submanifold of \(L M\).  However, although in finite dimensions
geodesic completeness is equivalent to a lot of things, that no longer
holds for loop spaces.  In particular, although \(\exp : T_x M \to M\)
is surjective, it need not be the case that \(\exp : T_\alpha L M \to
L M\) is surjective.  For example, take the sphere, \(S^2\), and the
exponential map based at the constant loop at the south pole.
Consider a loop which is a great circle through the south pole.  When
we try to lift this to the tangent space at the south pole, we find
that it must lift to a segment of a straight line between preimages of
the south pole, as it is a geodesic segment.  This cannot be made into
a loop and so there is no lift.  Thus there is no geodesic between the
constant loop at the south pole and any great circle through this
point.

\subsection{More Fun with Based Loops}

We conclude with a brief introduction to the topic of
\emph{polarisations}.  For simplicity, we shall work with complex
vector bundles.

The space \(L \C\) has a very simple description using Fourier
analysis.  It is an appropriate completion of the space of Laurent
polynomials, \(\C[z,z^{-1}]\).  In particular it decomposes as two
pieces: \(L_- \C \oplus L_+ \C\) according to the powers of \(z\).  We
assign the constant loops, corresponding to \(z^0\), to \(L_+ \C\).
These have more stylish descriptions as the space of loops that extend
holomorphically over an inner or outer disc (modulo the assignment of
the constant loops).

One might ask whether this structure is preserved on a loop space.
That is, given a complex vector bundle \(E \to M\), is there a similar
fibrewise splitting of \(L E\)?  If one restricts \(L E\) to the
constant loops then this does exist, but over the whole of \(L M\)
then it does not except in very special circumstances (see~\cite{rcas}
and~\cite{as2}).

However, it almost works.  When moving from one chart to another one
finds that the projections \(L_\pm \C^n \to \widetilde{L_\pm} \C^n\)
are Fredholm and \(L_\pm \C^n \to \widetilde{L_\mp} \C^n\) are
compact.  This means that, morally, one is only shifting a finite
dimensional amount from one side to the other.

Such a structure is called a \emph{polarisation}.  The definitive
reference is~\cite{apgs}.  The theory of polarisations is intimately
connected with that of representations of loop groups which is why it
is of particular interest to students of loop spaces.  We shall just
mention a few highlights here.

The reason for the title of this section is that the based loop space,
\(\Omega U_n\) is closely linked to polarisations.  Given a complex
vector bundle \(E \to M\), one can consider the bundle over \(L M\)
the points of which are the splittings of \(L E\) which are equivalent
to the canonical polarisation.  These are only fibrewise splittings so
always exist.  A global section of this bundle would define a global
splitting of \(L E\) which, by work of~\cite{rcas} and~\cite{as2},
would mean that \(L E \to L M\) was ``almost'' trivial.  If one
imposes a little extra structure on the splittings, namely that they
are orthogonal and behave well under the natural \(L \C\)-action, then
this bundle is very easy to identify: let \(Q \to M\) be the principal
\(U_n\)-bundle associated to \(E\).  The group \(L U_n\) acts on
\(\Omega U_n\) via \(\gamma \cdot \beta = \gamma \beta
\gamma(0)^{-1}\).  The bundle of ``nice'' polarisations is:
\[
L Q \times_{L U_n} \Omega U_n.
\]
The bundle of \emph{all} polarisations is homotopy equivalent to \(L Q
\times_{L U_n} \Omega U\).

One can give an alternative interpretation of this bundle.  A point of
\(L Q\) consists of a trivialisation of \(\alpha^* E\) for some
\(\alpha \in L M\).  That is, it is a fibrewise isomorphism \(\alpha^*
E \cong S^1 \times \C^k\).  A point of \(L Q \times_{L U_n} \Omega
U_n\) is also a trivialisation of \(\alpha^* E\) but not to a
``standard'' reference space such as \(\C^k\).  Rather, it is a
fibrewise isomorphism \(\alpha^* E \cong S^1 \times E_{\alpha(0)}\).
Thus a global section of this bundle defines an isomorphism \(L E
\cong e_0^* E \otimes L \C\) which is what is meant by a bundle being
``almost'' trivial.

Another highlight of the topic of polarisations is that although the
subspaces \(L_+ E\) and \(L_- E\) are not well-defined, the exterior
algebra \(\exterior^\bullet (L_+ E)^* \widetilde{\otimes}
\exterior^\bullet L_- E\) is well-defined%
\footnote{%
  The tilde on the tensor product is to denote an appropriate
  completion.
}.
The grading here is slightly odd in that the \(L_- E\) part is
negatively graded.  Thus the degree of \(\exterior^{p} (L_+ E)^*
\widetilde{\otimes} \exterior^q L_- E\) is \(p - q\).  This is known
as the \emph{semi-infinite exterior power} of \(L E\).  This is
because in finite dimensions a choice of isomorphism \(\exterior^{\dim
W} W \cong \C\) defines an isomorphism:
\[
\exterior^{\bullet + \dim W} V^* \cong \exterior^\bullet W^0 \otimes
\exterior^\bullet W
\]
for \(W \subseteq V\) where \(W^0\) is the annihilator of \(W\) in
\(W^*\) and \(\exterior^\bullet W\) is negatively graded.  Hence as
\(L E\) is infinite dimensional and \(L_- E\) is of half the dimension
of \(L E\), we get:
\[
\exterior^{\bullet + \infty/2} (L E)^* \cong \exterior^\bullet (L_+
E)^* \otimes \exterior^\bullet L_- E.
\]

One simple reason for wanting to take account of such structure can be
seen from the idea of the signature of a manifold.  This turns out to
depend on the middle dimension cohomology.  If one wanted to
generalise this to loop spaces, one would need a notion of ``middle
dimension'' cohomology, i.e.~of \emph{semi-infinite} cohomology.  A
deeper reason is that this is part of the link between polarisations
and representations of loop groups, about which I shan't go into more
detail except to say that it is closely linked to the theory of spin
and spin bundles on loop spaces.  A good place to start reading is the
book~\cite{apgs} and a good place to end reading is the recent result
of Freed, Hopkins, and Teleman~\cite{math.AT/0312155}.  En route, I
would take in~\cite{rppr} and my preprint on the construction of the
Dirac operator in~\cite{math.DG/0505077}.


\newpage

\begin{thebibliography}{Bom67}

\bibitem[Bom67]{jb3}
Jan Boman.
\newblock Differentiability of a function and of its compositions with
  functions of one variable.
\newblock {\em Math. Scand.}, 20:249--268, 1967.

\bibitem[CG04]{rcvg}
Ralph~L. Cohen and V{\'e}ronique Godin.
\newblock A polarized view of string topology.
\newblock In {\em Topology, geometry and quantum field theory}, volume 308 of
  {\em London Math. Soc. Lecture Note Ser.}, pages 127--154. Cambridge Univ.
  Press, Cambridge, 2004.

\bibitem[CJ02]{rcjj}
Ralph~L. Cohen and John D.~S. Jones.
\newblock A homotopy theoretic realization of string topology.
\newblock {\em Math. Ann.}, 324(4):773--798, 2002.

\bibitem[CS04]{rcas}
Ralph~L. Cohen and Andrew Stacey.
\newblock Fourier decompositions of loop bundles.
\newblock In {\em Homotopy theory: relations with algebraic geometry, group
  cohomology, and algebraic $K$-theory}, volume 346 of {\em Contemp. Math.},
  pages 85--95. Amer. Math. Soc., Providence, RI, 2004.

\bibitem[FHT]{math.AT/0312155}
Daniel~S. Freed, Michael~J. Hopkins, and Constantin Teleman.
\newblock {Twisted K-theory and Loop Group Representations},
  arXiv:math.AT/0312155.

\bibitem[Jon87]{jj}
John D.~S. Jones.
\newblock Cyclic homology and equivariant homology.
\newblock {\em Invent. Math.}, 87(2):403--423, 1987.

\bibitem[KM97]{akpm}
Andreas Kriegl and Peter~W. Michor.
\newblock {\em The convenient setting of global analysis}, volume~53 of {\em
  Mathematical Surveys and Monographs}.
\newblock American Mathematical Society, Providence, RI, 1997.

\bibitem[LM89]{hlmm}
H.~Blaine Lawson, Jr. and Marie-Louise Michelsohn.
\newblock {\em Spin geometry}, volume~38 of {\em Princeton Mathematical
  Series}.
\newblock Princeton University Press, Princeton, NJ, 1989.

\bibitem[PR94]{rppr}
R.~J. Plymen and P.~L. Robinson.
\newblock {\em Spinors in {H}ilbert space}, volume 114 of {\em Cambridge Tracts
  in Mathematics}.
\newblock Cambridge University Press, Cambridge, 1994.

\bibitem[PS86]{apgs}
Andrew Pressley and Graeme Segal.
\newblock {\em Loop groups}.
\newblock Oxford Mathematical Monographs. The Clarendon Press Oxford University
  Press, New York, 1986.
\newblock Oxford Science Publications.

\bibitem[Sta]{math.DG/0505077}
Andrew Stacey.
\newblock {The Geometry of the Loop Space and a Construction of a Dirac
  Operator}, arXiv:math.DG/0505077.

\bibitem[Sta05a]{as2}
Andrew Stacey.
\newblock Finite dimensional subbundles of loop bundles.
\newblock {\em Pacific J. Math.}, 219(1):187--199, 2005.

\bibitem[Sta05b]{as3}
Andrew Stacey.
\newblock The truncated {W}itten genus.
\newblock {\em Math. Z.}, 249(3):581--595, 2005.

\end{thebibliography}
\end{document}